\documentclass[11pt]{article}
\usepackage[top=2.54cm, bottom=2.54cm, left=2.5cm, right=2.5cm]{geometry}

\usepackage[T1]{fontenc}
\usepackage[utf8]{inputenc}
\usepackage{lmodern}
\usepackage{microtype}

\usepackage[tbtags]{mathtools}
\usepackage{amssymb, latexsym}
\usepackage{amsthm}
\usepackage{mathrsfs}

\usepackage{enumitem}
\usepackage{aliascnt}

\usepackage{xcolor}
\usepackage[numbers,sort&compress]{natbib}
\usepackage[colorlinks=true, linkcolor=blue, citecolor=blue, urlcolor=magenta]{hyperref}
\usepackage[nameinlink,capitalize]{cleveref}

\usepackage{lineno}

\usepackage{graphicx}
\usepackage{appendix}
\usepackage{fancyhdr}

\usepackage{float}
\usepackage{pict2e}
\usepackage{tikz}
\usepackage{tikz-3dplot}
\usepackage{tkz-graph}
\usetikzlibrary{arrows,shapes,chains}
\usepackage {verbatim}
\usepackage{scalerel}

\allowdisplaybreaks[4]
\numberwithin{equation}{section}

\theoremstyle{plain}
\newtheorem{theorem}{Theorem}[section]

\newaliascnt{lemma}{theorem}
\newtheorem{lemma}[lemma]{Lemma}
\aliascntresetthe{lemma}

\newaliascnt{corollary}{theorem}
\newtheorem{corollary}[corollary]{Corollary}
\aliascntresetthe{corollary}

\theoremstyle{remark}

\newaliascnt{remark}{theorem}

\aliascntresetthe{remark}

\theoremstyle{definition}

\newaliascnt{definition}{theorem}
\newtheorem{definition}[definition]{Definition}
\aliascntresetthe{definition}

\newaliascnt{conjecture}{theorem}

\aliascntresetthe{conjecture}

\newaliascnt{problem}{theorem}
\newtheorem{problem}[problem]{Problem}
\aliascntresetthe{problem}

\newaliascnt{proposition}{theorem}

\aliascntresetthe{proposition}

\newaliascnt{claim}{theorem}

\aliascntresetthe{claim}

\newaliascnt{observation}{theorem}

\aliascntresetthe{observation}

\newaliascnt{fact}{theorem}

\aliascntresetthe{fact}

\crefname{theorem}{Theorem}{Theorems}
\crefname{lemma}{Lemma}{Lemmas}
\crefname{corollary}{Corollary}{Corollaries}

\crefname{remark}{Remark}{Remarks}
\crefname{problem}{Problem}{Problems}
\crefname{conjecture}{Conjecture}{Conjectures}

\crefname{definition}{Definition}{Definitions}
\crefname{proposition}{Proposition}{Propositions}
\crefname{claim}{Claim}{Claims}
\crefname{observation}{Observation}{Observations}
\crefname{fact}{Fact}{Facts}

\usepackage{nomencl}
\makenomenclature

\newcommand{\N}{\mathbb N}
\newcommand{\eps}{\varepsilon}
\newcommand{\A}{\mathcal A}
\newcommand{\B}{\mathcal B}
\newcommand{\C}{\mathcal C}
\newcommand{\Pcal}{\mathcal P}

\newcommand{\deltacd}{\delta_{\ee}}

\newcommand{\bb}{\boldsymbol{b}}

\newsavebox\vvvbox
\savebox\vvvbox{\tikz{
		\draw[black,fill=black] (90:1) circle (.35);
		\draw[black,fill=black] (210:1) circle (.35);
		\draw[black,fill=black] (330:1) circle (.35);
		\draw[opacity=0] (0:1.2) circle (0.1);
}}
\newcommand{\vvv}{\mathord{\scaleobj{1.2}{\scalerel*{\usebox{\vvvbox}}{x}}}}
\newcommand{\pivvv}{\pi_{\vvv}}

\newsavebox\evbox
\savebox\evbox{\tikz{
		\draw[black,fill=black] (90:1) circle (.35);
		\draw[black,fill=black] (210:1) circle (.35);
		\draw[black,fill=black] (330:1) circle (.35);
		\draw[black,line width=0.28cm ] (210:1) -- (330:1);
		\draw[opacity=0] (0:1.2) circle (0.1);
}}
\newcommand{\ev}{\mathord{\scaleobj{1.2}{\scalerel*{\usebox{\evbox}}{x}}}}	
\newcommand{\piev}{\pi_{\ev}}

\newsavebox\eebox
\savebox\eebox{\tikz{
		\draw[black,fill=black] (90:1) circle (.35);
		\draw[black,fill=black] (210:1) circle (.35);
		\draw[black,fill=black] (330:1) circle (.35);
		\draw[black,line width=0.28cm ] (90:1) -- (330:1);
		\draw[black,line width=0.28cm ] (90:1) -- (210:1);
		\draw[opacity=0] (0:1.2) circle (0.1);
}}
\newcommand{\ee}{\mathord{\scaleobj{1.2}{\scalerel*{\usebox{\eebox}}{x}}}}
\newcommand{\piee}{\pi_{\ee}}

\newcommand{\pieepal}{\pi_{\ee}^{\mathsf{pal}}}
\newcommand{\pievpal}{\pi_{\ev}^{\mathsf{pal}}}
\newcommand{\rev}{\operatorname{rev}}

\title{Vanishing Uniform Tur\'an Densities: Characterization and Jumps}

\author{Hao Lin \thanks{School of Mathematical Sciences, University of Science and Technology of China, Hefei, Anhui 230026, China. Email: \texttt{haolinz6@ustc.edu.cn}.} \and
Jie Ma \thanks{School of Mathematical Sciences, University of Science and Technology of China, Hefei, Anhui 230026, and Yau Mathematical Sciences Center, Tsinghua University, Beijing 100084, China. Research supported by National Key Research and Development Program of China 2023YFA1010201 and National Natural Science Foundation of China grant 12125106. Email: \texttt{jiema@ustc.edu.cn}.} \and
Wenling Zhou \thanks{School of Mathematics, Shandong University, Jinan, Shandong 250100, China. Research supported by the National Natural Science Foundation of China (No. 12401457),  Shandong Provincial Natural Science Foundation (No. ZR2024QA067), China Postdoctoral Science Foundation (No. 2024M761780), and Young Talent of Lifting engineering for Science and Technology in Shandong, China (No. SDAST2025QTA074). Email: \texttt{gracezhou@sdu.edu.cn}.
\\
\indent \textbf{Keywords:}
Uniform Tur\'an density, max-layered $3$-graphs, jump points.}
}

\date{}

\begin{document}
\maketitle

\begin{abstract}
Uniform Tur\'an density problems ask for the minimum density threshold that guarantees a copy of a fixed $3$-uniform hypergraph $F$ in every sufficiently large $3$-uniform hypergraph whose edges are uniformly distributed. This line of research was initiated by Erd\H{o}s and S\'os in the 1980s, who introduced the notion of \emph{vertex-set Tur\'an density} $\pivvv(F)$. Reiher, R\"odl and Schacht later defined two related concepts, namely, the \emph{vertex-pair Tur\'an density} $\piev(F)$ and the \emph{cherry Tur\'an density} $\piee(F)$, which satisfy the chain of inequalities
$\pivvv(F)\ge\piev(F)\ge\piee(F).$   
In this paper, we investigate the problem for their vanishing densities. 

First, we characterize the $3$-uniform hypergraphs $F$ for which
$\piee(F)=0$, thereby answering a question posed by Schacht. We show that these are precisely the so-called \emph{max-layered} $3$-uniform hypergraphs recently introduced by Lamaison.
Furthermore, we answer an analogous question of Lamaison for the \emph{minimum-codegree palette density} $\pieepal$. Together, these results imply that the following
are equivalent: $\piee(F)=0$, $\pieepal(F)=0$, and $F$ is a max-layered $3$-uniform hypergraph. 

Let $\star \in \{\vvv, \ev, \ee \}$. We say that zero is a \emph{jump} for $\pi_{\star}$ if there exists $c>0$ such that no $3$-uniform hypergraph $F$ satisfies
$\pi_{\star}(F)\in(0,c)$. 
A famous result of Reiher, R\"odl and Schacht asserts that zero is a jump for $\pivvv$.
In contrast, we prove that zero is not a jump for either $\piev$ or $\piee$. This also yields a negative answer to Schacht's well-ordering question for $\piev$.
Together with existing results, these findings complete the picture for all three notions of uniform Tur\'an density at zero.
\end{abstract}


\section{Introduction}
A central problem in Extremal Combinatorics, initiated by Tur\'an~\cite{Turan1941}, asks how dense a large graph or hypergraph can be while avoiding a fixed subgraph. For an integer $k\ge 2$, a \emph{$k$-uniform hypergraph}, or simply a \emph{$k$-graph}, is a pair $H=(V,E)$ with $E\subseteq \binom{V}{k}$. 
The \emph{Tur\'an number} $\operatorname{ex}(n,F)$ of a $k$-graph $F$ denotes the maximum number of edges in an $n$-vertex $k$-graph that is {\it $F$-free}, i.e., containing no copy of $F$ as a subgraph. 
The standard averaging argument shows that $\operatorname{ex}(n,F)/\binom nk$ is non-increasing in $n$, so we may legitimately define the \emph{Tur\'an density} of $F$ by
\[
\pi(F)=\lim_{n\to\infty}\frac{\operatorname{ex}(n,F)}{\binom nk}.
\]

The study of Tur\'an densities has been the focus of extensive research. 
For every $k\geq 2$, let $\Pi^{(k)}=\{\pi(F): F \text{ is a } k\text{-graph}\}$ consist of all Tur\'an densities of
single $k$-graphs. 
For graphs, Tur\'an densities are well understood. 
Mantel's theorem~\cite{mantel1907problem} and Tur\'an's theorem~\cite{Tu-graph} yield the exact Tur\'an numbers of complete graphs, while the celebrated Erd\H{o}s--Stone--Simonovits theorem~\cite{ErdosStone1946,ErdosSimonovits1966} determines the Tur\'an density $\pi(F)=1-\frac{1}{\chi(F)-1}$
for every graph $F$ with chromatic number $\chi(F)\ge 2$. 
The latter result implies that $\Pi^{(2)}=\{(r-1)/r: r\in \mathbb N^+ \}$.
For $k\ge 3$, however, it is notoriously difficult to understand the Tur\'an density $\pi(F)$ for general $k$-graphs $F$; see Keevash~\cite{Keevash11} for a survey.

Throughout the paper, we focus on $3$-graphs.
Even the Tur\'an densities of two of the simplest $3$-graphs,
the complete $3$-graph $K^{(3)}_4$ on four vertices and the
$3$-graph $K^{(3)-}_4$ obtained from it by deleting one edge,
remain unknown~\cite{k43-upbound,k43lowbound,k43-lowbound}.
Very little has been known about $\Pi^{(3)}$; see,
e.g.,~\cite{conlon2025hypergraphs}.

It is historically important to understand the behavior of Tur\'an densities at zero. There are two basic questions.
The first is to characterize hypergraphs $F$ for which $\pi(F)=0$. The second is to determine whether the positive Tur\'an densities are bounded away from zero or can be arbitrarily close to zero. More precisely, we say that zero is a \emph{jump} in $\Pi^{(3)}$ if there exists
$c>0$ such that  $\Pi^{(3)}\cap (0, c)=\emptyset$.
For the classical Tur\'an density of $3$-graphs, the first question
was settled by Erd\H{o}s~\cite{r-partite}, and the second follows
from his characterization.
A $3$-graph is \emph{$3$-partite} if its
vertex set can be partitioned into three classes such that every edge meets each class in exactly one vertex.

\begin{theorem}[Erd\H{o}s~\cite{r-partite}]
For every $3$-graph $F$, $\pi(F)=0$ if and only if $F$ is $3$-partite.
\end{theorem}

Note that if $F$ is not $3$-partite, then every complete $3$-partite $3$-graph is $F$-free.
Taking the balanced complete $3$-partite $3$-graph gives an
$F$-free $3$-graph of asymptotic edge density $2/9$. 
This implies that $\pi(F)>0$ in fact forces $\pi(F)\ge\frac{2}{9}.$
In particular, this shows that zero is a jump in $\Pi^{(3)}$. 

Throughout the paper, for distinct vertices $x,y,z$, we write $xy$ for $\{x,y\}$ and
$xyz$ for $\{x,y,z\}$. 
For a positive integer $\ell$, we denote by $[\ell]$ the set $\{1, \dots, \ell\}$.

\subsection{Uniform Tur\'an densities: background and their behavior at zero}

Many extremal constructions for classical hypergraph Tur\'an
problems have highly non-uniform edge distributions and often contain
large independent sets. This led Erd\H{o}s and S\'os~\cite{ErdosSos1982,Erdos-1990} in the 1980s to propose studying Tur\'an problems under the additional
requirement that the edges of the host $3$-graph be uniformly
distributed.
In this line of research, three natural notions of uniform Tur\'an
density for $3$-graphs have been introduced and have attracted
considerable attention: the \emph{vertex-set}, \emph{vertex-pair},
and \emph{cherry} Tur\'an densities.

We begin with the vertex-set Tur\'an density.
For
$d\in[0,1]$ and $\mu>0$, a $3$-graph $H=(V,E)$ is
\emph{$(d,\mu,\vvv)$-dense} if, for all $X,Y,Z\subseteq V$, the number
$e_{\vvv}(X,Y,Z)$ of ordered triples $(x,y,z)\in X\times Y\times Z$ with
$xyz\in E$ satisfies
\[
 e_{\vvv}(X,Y,Z)
 \ge d|X||Y||Z|-\mu|V|^3.
\]
The corresponding \emph{vertex-set Tur\'an density} of a fixed
$3$-graph $F$ is defined by
\[
\pivvv(F)=\sup\Bigl\{d\in[0,1]:
\begin{array}{l}
 \text{for every $\mu>0$ and $n_0\in\N$, there exists an $F$-free,}\\
 \text{$(d,\mu,\vvv)$-dense $3$-graph $H$ with $|V(H)|\ge n_0$}
\end{array}
\Bigr\}.
\]
With this notation,
Erd\H{o}s and S\'os asked to determine $\pivvv(K_4^{(3)})$ and
$\pivvv(K_4^{(3)-})$. The conjecture that $\pivvv(K^{(3)-}_4)=1/4$ was proved by Glebov, Kr\'al' and Volec~\cite{GlebovKralVolec2016}, and independently by Reiher, R\"odl and Schacht~\cite{ReiherRodlSchachtJEMS2018}. 
More recently, the conjecture $\pivvv(K_4^{(3)})=1/2$ was proved by Kielak, Kr\'al',
Lamaison, Liu, Shu and Wu~\cite{k43-1}, and independently by Buci\'c~\cite{k43-2}, confirming the
optimality of R\"odl's construction~\cite{Rodl1986}.
Besides, recent work has also determined several new exact values for
$\pivvv$ including the minimum
positive value~\cite{1/27}, tight cycles~\cite{cycle-uniform}, large
stars~\cite{LamaisonWu2024,LinZhouStars2025} and families with density
$1/4$~\cite{chen2022beyond,Hypergraph-llwzhou23} or
$8/27$~\cite{8/27}. Further structural and
classification results can be found in~\cite{ Lamaison2024, KKLT-25, LWZZ2025,King-piga-SS-25}.

The $\vvv$-dense condition provides a natural perspective to avoid extremal constructions with large independent sets. 
In their treatment of hypergraph Tur\'an problems,
Reiher, R\"odl and Schacht~\cite{ReiherRodlSchachtk43,RRS-cherry}
introduced two stronger notions of uniform density.
For $d\in[0,1]$ and $\mu>0$, a $3$-graph $H=(V,E)$ is
\emph{$(d,\mu,\ev)$-dense} if, for any subset $X\subseteq V$ and any set of ordered pairs $P\subseteq V\times V$, the number $e_{\ev}(X,P)$ of elements $(x,(y,z))\in X\times P$
with $xyz\in E$ satisfies
\[
 e_{\ev}(X,P)
 \ge d|X||P|-\mu|V|^3.
\]
The corresponding \emph{vertex-pair Tur\'an density} of a fixed
$3$-graph $F$ is defined by
\[
\piev(F)=\sup\Bigl\{d\in[0,1]:
\begin{array}{l}
 \text{for every $\mu>0$ and $n_0\in\N$, there exists an $F$-free,}\\
\text{$(d,\mu,\ev)$-dense $3$-graph $H$ with $|V(H)|\ge n_0$}
\end{array}
\Bigr\}.
\]
The vertex-pair condition can be strengthened once more by allowing
the two ordered pairs $(x,y)$ and $(y,z)$, which share the middle
vertex $y$. More precisely, a $3$-graph $H=(V,E)$ is called \emph{$(d,\mu,\ee)$-dense} if for every $P, Q\subseteq V\times V$, the number $e_{\ee}(P,Q)$ of members $((x,y),(y,z))\in P\times Q$
with $xyz\in E$ satisfies
\[
e_{\ee}(P,Q) \ge d|\mathcal K_{\ee}(P,Q)|-\mu |V|^3,
\]
where $\mathcal K_{\ee}(P,Q)=\{((x,y),(y,z))\in P\times Q\}$.
The corresponding \emph{cherry Tur\'an density} of a given
$3$-graph $F$ is defined by
\[
\piee(F)=\sup\Bigl\{d\in[0,1]:
\begin{array}{l}
 \text{for every $\mu>0$ and $n_0\in\N$, there exists an $F$-free,}\\
 \text{$(d,\mu,\ee)$-dense $3$-graph $H$ with $|V(H)|\ge n_0$}
\end{array}
\Bigr\}.
\]
By the above definitions, every $(d,\mu,\ee)$-dense $3$-graph is also $(d,\mu,\ev)$-dense, and every $(d,\mu,\ev)$-dense $3$-graph is also $(d,\mu,\vvv)$-dense. Therefore, for every $3$-graph
$F$,
\begin{equation}
\pivvv(F)\ge\piev(F)\ge\piee(F).   
\end{equation}

For the tetrahedron $K_4^{(3)}$, all three uniform Tur\'an
densities considered here are now known. The result \cite{k43-1,k43-2} mentioned above, together with the main theorem in~\cite{ReiherRodlSchachtk43}, gives
$\pivvv(K_4^{(3)})=\piev(K_4^{(3)})=1/2$, whereas
$\piee(K_4^{(3)})=0$ was proved in~\cite{RRS-cherry}.
More generally, the results
in~\cite{RRS-cherry,BergerPigaReiherRodlSchachtK5} determine
$\piee(K_t^{(3)})$ for every $4\le t\le16$, except for $t\in \{9,10\}$.

We now return to the two questions concerning the vanishing densities introduced above. 
Reiher, R\"odl and Schacht~\cite{ReiherRodlSchachtJLMS2018} gave a characterization of $3$-graphs with vertex-set
Tur\'an density being zero.
To state it, let $\partial F:=\{uv\in\binom{V(F)}2: uv\subseteq e\text{ for some }e\in E(F)\}$ be the \emph{shadow} of a $3$-graph $F$.

\begin{theorem}[Reiher, R\"odl and
Schacht~\cite{ReiherRodlSchachtJLMS2018}]\label{thm:dot-zero}
For every $3$-graph $F$, $\pivvv(F)=0$ if and only if there
 exist a linear order $\prec$ on $V(F)$, and a $3$-coloring $\varphi:\partial F\to \{\mathrm{red},\mathrm{blue},\mathrm{green}\}$
 such that every edge $uvw\in E(F)$ with $u\prec v \prec w$ satisfies
 \[
\varphi(uv)=\mathrm{red},\qquad
\varphi(uw)=\mathrm{blue},\qquad
\varphi(vw)=\mathrm{green}.
 \]
\end{theorem}
If a $3$-graph $F$ satisfies $\pivvv(F)>0$, then its vertex-set
Tur\'an density is at least $1/27$ because of the following construction from~\cite{ReiherRodlSchachtJLMS2018}. 
Consider a random $n$-vertex $3$-graph
$H_n$ as follows: let $V(H_n)=[n]$; color each pair of $[n]$ red, blue, or
green independently and uniformly at random, and include $ijk$ with $1\le i < j < k \le n$ as an edge of $H_n$ if $ij, ik, jk$ have colors red, blue, and green, respectively.
By~\cref{thm:dot-zero}, it is straightforward to see that $H_n$ is
$F$-free if $\pivvv(F)> 0$. 
On the other hand, since each triple becomes an edge with
probability $1/27$, a standard concentration argument shows that, with high probability,
$H_n$ is $(1/27,o(1),\vvv)$-dense. Thus, zero is a jump for $\pivvv$.
In fact, $1/27$ is the smallest positive density attained by $\pivvv$ (see Garbe, Kr\'al' and Lamaison~\cite{1/27}). 

Thus, for $\pivvv$, both the zero‑density class and the behavior immediately above zero are completely characterized. It is natural to ask for analogous results for the
two stronger densities $\piev$ and $\piee$. In his 2022 ICM report,
Schacht~\cite{SchachtICM} 
posed three problems concerning their
behavior at zero. One of these asked for a useful characterization
of the $3$-graphs $F$ with
$\piev(F)=0$; see~\cite[Problem~4.4]{SchachtICM}. This question was
subsequently answered in the breakthrough work of Lamaison~\cite{Lamaison2024}.

To state Lamaison's characterization, we first
introduce the following notation.

\begin{definition}[Edge-patterns]
	\label{def:intro-edge-patterns}
	Let $F$ be a $3$-graph, let $\prec$ be an order on $V(F)$, and let
	$\lambda:\partial F\to\N$.  For an edge $uvw\in E(F)$ with
	$u\prec v\prec w$, the ordered triple $\bigl(\lambda(uv),\lambda(uw),\lambda(vw)\bigr)$
	is called the \emph{edge-pattern} of $uvw$ with respect to
$(\prec,\lambda)$. The \emph{edge-pattern set} of $F$ with respect to
$(\prec,\lambda)$ is denoted by
\begin{equation}
		\label{edge-pattern-set}
		\B_{\lambda,\prec}
		:=
		\left\{
		\bigl(\lambda(uv),\lambda(uw),\lambda(vw)\bigr):
		uvw\in E(F),\ u\prec v\prec w
		\right\}.
	\end{equation}
\end{definition}

A set $B\subseteq\N^3$ is \emph{min-degenerate} if every member of $B$ has three distinct entries and, whenever two distinct members of $B$ share an entry, that entry is the minimum entry of at least one of them. A $3$-graph $F$ is \emph{min-layered} if there exist an order $\prec$ on $V(F)$ and a map $\lambda:\partial F\to\N$ such that $\B_{\lambda,\prec}$ is min-degenerate.
Lamaison~\cite{Lamaison2024} proved the following characterization.

\begin{theorem}[Lamaison~\cite{Lamaison2024}]\label{thm:zero-dot-edge}
For every $3$-graph $F$, $\piev(F)=0$ if and only if $F$ is min-layered.
\end{theorem}

Let $\Pi^{(3)}_{\ev}
= \{\piev(F):F\text{ is a }3\text{-graph}\}$.
Although \cref{thm:zero-dot-edge} gives a zero‑density characterization for $\piev$, it remained open whether zero is a jump in $\Pi^{(3)}_{\ev}$.
Another problem of Schacht~\cite{SchachtICM} asked a more general question for the set $\Pi^{(3)}_{\ev}$ as follows. 
We say a subset $S \subseteq\mathbb R$ is \emph{well-ordered} under the usual order if every nonempty subset of $S$ has a least element.

\begin{problem}[\cite{SchachtICM}, Problem~4.3(iii)]
\label{prob:schacht-1}
Show that $\Pi^{(3)}_{\ev}$
	is well-ordered as a subset of the reals.
\end{problem}

\noindent A positive answer to this problem would, in particular,
imply that zero is a jump in $\Pi^{(3)}_{\ev}$.

The last problem of Schacht~\cite{SchachtICM} mentioned above asked for a characterization of the zero-density class of the cherry Tur\'an density.

\begin{problem}[\cite{SchachtICM}, Problem~5.3(i)]
	\label{prob:schacht-3}
	Find a useful characterization of the $3$-graphs $F$ with
	$\piee(F)=0$.
\end{problem}

Lamaison~\cite{Lamaison2024} had also asked the analogous zero-density question for the
minimum-codegree palette density (see the precise definition of $\pieepal(F)$ we give shortly).

These questions stated above provide the driving motivation for our investigation of vanishing uniform Tur\'an densities.
We will introduce our main results in the next subsection, which, together with the existing results, give a complete picture for all three notions of uniform Tur\'an density at zero.

\subsection{Our results}

Our first result characterizes the zero-density class of the cherry
Tur\'an density and answers \cref{prob:schacht-3}.  It is stated in terms of
max-layeredness, a notion introduced by Lamaison~\cite{Lamaison2024}.  A set
$B\subseteq\N^3$ is \emph{max-degenerate} if the maximum entry of every
member of $B$ occurs in a unique position and distinct members of $B$ have
distinct maximum entries.

\begin{definition}[Max-layered $3$-graphs~\cite{Lamaison2024}]
	\label{def:intro-max-layered}
	A $3$-graph $F$ is \emph{max-layered} if there exist an order $\prec$
	on $V(F)$ and a map $\lambda:\partial F\to\N$ such that the edge-pattern
	set $\B_{\lambda,\prec}$ is max-degenerate.  In this case,
	$(\prec,\lambda)$ is called a \emph{max-layered witness} for $F$.
\end{definition}

\begin{theorem}\label{thm:main}
For every $3$-graph $F$, $\piee(F)=0$ if and only if $F$ is max-layered.
\end{theorem}

Our argument also settles a question of Lamaison~\cite{Lamaison2024}. A \emph{palette} is a pair $\mathscr P=(\mathscr C,\mathscr A)$, where
$\mathscr C$ is a finite nonempty set of colors and
$\mathscr A\subseteq\mathscr C^3$ is a set of admissible ordered color
triples.  A $3$-graph $F$ \emph{admits} $\mathscr P$ if there exist an order
$\prec$ on $V(F)$ and a coloring $\varphi:\partial F\to\mathscr C$ such that
\[
\bigl(\varphi(uv),\varphi(uw),\varphi(vw)\bigr)\in\mathscr A
\]
for every edge $uvw\in E(F)$ with $u\prec v\prec w$.  We call
$(\prec,\varphi)$ an \emph{admission witness} for $F$ and $\mathscr P$. 
In fact, if we view $\mathscr{C}$ as $[|\mathscr{C}|]$, then $F$ \emph{admits} $\mathscr{P}$ if and only if the edge‑pattern $\mathcal{B}_{\varphi, \prec}\subseteq \mathscr{A}$.

For $a,b\in\mathscr C$, define $d_{12}(a,b):=\bigl|\{c\in\mathscr C:(a,b,c)\in\mathscr A\}\bigr|,
 d_{13}(a,b):=\bigl|\{c\in\mathscr C:(a,c,b)\in\mathscr A\}\bigr|,
 d_{23}(a,b):=\bigl|\{c\in\mathscr C:(c,a,b)\in\mathscr A\}\bigr|$.
The \emph{minimum codegree} of $\mathscr P$ is
\begin{equation}
   \deltacd(\mathscr P)
:=
\frac{1}{|\mathscr C|}
\min_{a,b\in\mathscr C}
\min\{d_{12}(a,b),d_{13}(a,b),d_{23}(a,b)\}. 
\end{equation}
The corresponding \emph{minimum-codegree palette density} is
\[
\pieepal(F)
:=
\sup\bigl\{\deltacd(\mathscr P):
\mathscr P\text{ is a palette not admitted by }F\bigr\}.
\]
Lamaison~\cite{Lamaison2024} asked whether $\pieepal(F)=0$ is equivalent to that $F$ is max-layered. We answer this question affirmatively.

\begin{theorem}\label{thm:palette-zero}
	For every $3$-graph $F$, $\pieepal(F)=0$ if and only if $F$ is max-layered.
\end{theorem}

Together, \cref{thm:main,thm:palette-zero} show that, for every
$3$-graph $F$, the following are equivalent:
\begin{equation}\label{eq:zero equivalent}
\piee(F)=0
\quad\Longleftrightarrow\quad
F\text{ is max-layered}
\quad\Longleftrightarrow\quad
\pieepal(F)=0.   
\end{equation}

We next turn to the densities immediately above zero.  In contrast to the
classical and vertex-set Tur\'an densities, the positive values of both
$\piev$ and $\piee$ can be arbitrarily close to zero.

\begin{theorem}\label{thm:zero-not-jump}
For every $\eps>0$, there exists a $3$-graph $F_\eps$ such that
$0<\piee(F_\eps)\le\piev(F_\eps)\le\eps$.
\end{theorem}

Thus, zero is not a jump for either $\piev$ or $\piee$.
For $\star \in \{\ev, \ee \}$, let $\Pi^{(3)}_{\star}=
\{\pi_{\star}(F):F\text{ is a }3\text{-graph}\}$. The above theorem also gives the following consequence.

\begin{corollary}\label{cor:not-well-ordered}
The sets $\Pi^{(3)}_{\ev}$ and $\Pi^{(3)}_{\ee}$
are not well-ordered subsets of the reals.
\end{corollary}

Indeed, \cref{thm:zero-not-jump} gives, for every $\eps>0$, a positive element of each set below $\eps$. Hence neither set has a least positive element. In particular, \cref{cor:not-well-ordered} gives a negative answer to
\cref{prob:schacht-1}.

Together with the existing results, our theorems determine the
zero-density classes (i.e., the families of $3$-graphs with vanishing density) and the jump behavior at zero for the three uniform Tur\'an densities considered here. The zero-density classes of $\pivvv$,
$\piev$, and $\piee$ are characterized, respectively, by the ordered $3$-coloring condition, min-layeredness, and max-layeredness.  Moreover, zero is a jump for $\pivvv$; in contrast, zero is not a jump for either $\piev$ or $\piee$. 

\subsection{Proof overview and organization}
We first present the proof overview of \cref{thm:main,thm:palette-zero},
then present the proof overview for \cref{thm:zero-not-jump}, and finally describe the organization of the paper.

\paragraph{The zero-density characterizations.}
Before we sketch the proofs of \cref{thm:main,thm:palette-zero}, we would like to point out that the two zero-density characterizations are closely related. We first prove that every max-layered
$3$-graph has zero density in both settings. For the
minimum-codegree palette density in \cref{thm:palette-zero}, this
follows from a direct greedy coloring of the layers. For the cherry
Tur\'an density in \cref{thm:main}, the same layer-by-layer idea is
carried out in the reduced $3$-graph formulation given by
\cref{thm:reiher-reduced}, which is the content of
\cref{lem:max-layered-zero}. 
The converse implications are proved together in
\cref{lem:nonmax-positive}. If $F$ is not
max-layered, we construct a palette of positive minimum codegree that
$F$ does not admit. This immediately gives
$\pieepal(F)>0$, and \cref{prop:palette-lower-bound} then yields
$\piee(F)>0$ as well. 
Together, these arguments establish both zero-density characterizations.
We now describe these proof overview in the following.

Suppose first that $F$ is max-layered. A palette has positive minimum codegree if,
once colors have been fixed in any two positions of an ordered color
triple, there is an admissible choice for the third. We give all pairs in the same layer the same color and
consider the layers in increasing order. When a new layer is reached,
max-layeredness ensures that there is at most one edge-pattern whose
maximum entry is that layer. 
The other two colors in this pattern have already been fixed, so positive minimum codegree provides a suitable color for the new layer.
Continuing in this way shows that $F$ admits every palette
of positive minimum codegree, proving the  corresponding direction of
\cref{thm:palette-zero}.

For the cherry Tur\'an density, the same layer-by-layer
strategy is carried out in reduced $3$-graphs
(\cref{def:reduced-3-graph}). Instead of assigning a color to each
layer, we choose representatives from the vertex classes of the
reduced $3$-graph. The reduced density condition provides many
possible choices, and a selection lemma is used to coordinate these
choices across all relevant triples. Together with
\cref{thm:reiher-reduced}, this proves the corresponding direction of
\cref{thm:main}; see \cref{sec:reduced-image} for more details.

It remains to show that, when $F$ is not max-layered, both densities
are positive. We construct a palette of positive minimum codegree that
$F$ does not admit. For each of the three positions, we prescribe the
missing color as a function of the other two, and declare a triple
admissible whenever it satisfies at least one of these three rules.
This automatically gives positive minimum codegree.

If $F$ admitted such a palette, an admission witness would give a
vertex order and a coloring of the pairs of $F$. Choosing, for each
edge, one of the three rules satisfied by its color triple would then
produce a finite system of equations in these pair colors. For each
fixed vertex order and choice of rules, the main lemma \cref{lem:functional-obstruction}
uses the failure of max-layeredness to choose the functions so that
the corresponding system has no solution. Since there are only
finitely many such choices, we combine the resulting constructions
coordinatewise to obtain a single palette that $F$ does not admit.
This construction is given in \cref{sec:palette-construction}.
It proves $\pieepal(F)>0$, and
\cref{prop:palette-lower-bound} then gives $\piee(F)>0$ as well; see \cref{sec:palette-construction} for more details.

\paragraph{The non-jump theorem.}
The proof of \cref{thm:zero-not-jump} also uses palettes, but in a
different way. Fix $\eta>0$. We first construct a universal palette
$\mathscr U_\eta$. A standard probability method shows that every palette whose
vertex-pair minimum degree is greater than $\eta$ receives a
homomorphism from a palette on a bounded number of colors, where a
palette homomorphism is a color map preserving admissible triples.
Taking the product of all such bounded palettes gives
$\mathscr U_\eta$, which has positive vertex-pair minimum degree and
maps homomorphically into every palette above the threshold $\eta$;
see \cref{lem:bounded-palette-reduction} and
\eqref{eq:universal-palette}. Consequently, if a $3$-graph $F$ admits
$\mathscr U_\eta$, then it admits every palette with vertex-pair
minimum degree greater than $\eta$. By
\cref{thm:lamaison-ev-palette}, this implies
$\piev(F)\le\eta$.

We then seek a positive lower bound for the cherry Tur\'an density of
the same $3$-graph. For this purpose, we construct the affine palette
$\mathscr Q_p$ from \cref{def:affine-palette}, whose minimum codegree
is $1/p$. By \cref{lem:affine-separator}, the prime $p$ can be chosen
so that $\mathscr U_\eta$ has no homomorphism to either
$\mathscr Q_p$ or its reverse. The palette separation theorem
\cref{thm:palette-separation} then gives a $3$-graph $F$ which admits
$\mathscr U_\eta$ but does not admit $\mathscr Q_p$. The first property
gives $\piev(F)\le\eta$, whereas the second and
\cref{prop:palette-lower-bound} give $\piee(F)\ge1/p>0$. 
Since $\eta$ can be chosen arbitrarily small, this proves the
non-jump theorem.

\paragraph{Organization.}
The remainder of the paper is organized as follows.  In
\cref{sec:preliminary}, we recall reduced $3$-graphs for the cherry
Tur\'an density.  In \cref{sec:main-proof}, we state the two key
lemmas (\cref{lem:max-layered-zero,lem:nonmax-positive}) and deduce \cref{thm:main,thm:palette-zero}.  The two key lemmas are proved in
\cref{sec:reduced-image,sec:palette-construction}, respectively.  In
\cref{sec:not-jump}, we prove \cref{thm:zero-not-jump} and obtain \cref{cor:not-well-ordered}.  Finally, in \cref{sec:conclusion},  we give some new explicit examples of $3$-graphs $F$ with $\piev(F)\neq \piee(F)$  and discuss further questions.

\section{Preliminaries}\label{sec:preliminary}
In this section, we recall reduced $3$-graphs and palette constructions for the cherry
Tur\'an density.

\subsection{Reduced \texorpdfstring{$3$}{3}-graphs}
\label{sec:cherry-reduced}

In this subsection, we introduce the concept of reduced $3$-graphs, which encapsulates hypergraph
regularity arguments about the uniform Tur\'an density of $3$-graphs. 
We follow some notions introduced in Reiher's survey~\cite[Section 3]{Reiher2020}.
It should be noted that the concept of reduced $3$-graphs has also been mentioned in other related works, e.g.,~\cite{ReiherRodlSchachtJEMS2018,cycle-uniform, 1/27,chen2022beyond, Hypergraph-llwzhou23}. 

\begin{definition}[Reduced $3$-graphs]\label{def:reduced-3-graph}
Let $I\subseteq\N$ be a finite index set.  For each
$ij\in\binom{I}{2}$, let $\Pcal^{ij}$ be a nonempty vertex set, and
assume that these sets are pairwise disjoint.  For every $ijk\in\binom{I}{3}$, let $\A^{ijk}$ be a $3$-partite $3$-graph with vertex partition $\{\Pcal^{ij},\Pcal^{ik}, \Pcal^{jk}\}$. Then the $\binom{|I|}{2}$-partite $3$-graph $\A$ with
\[
V(\A)=\bigcup_{ij\in\binom{I}{2}}\Pcal^{ij}
\qquad\text{and}\qquad
E(\A)=\bigcup_{ijk\in\binom{I}{3}}E(\A^{ijk})
\]
is called an \emph{$I$-reduced $3$-graph}.  The sets
$\Pcal^{ij}$ are its \emph{vertex classes}, and the $3$-graphs
$\A^{ijk}$ are its \emph{constituents}.  
\end{definition}

Throughout, the superscripts of vertex classes and constituents are
unordered: $\Pcal^{ij}=\Pcal^{ji}$, and $\A^{ijk}$ depends only on
the set $ijk$. Whenever an order is needed, the indices are
written increasingly.

\begin{definition}[Reduced $\ee$-density]\label{def:reduced-ee-density}
	Let $d\in[0,1]$ and let $\A$ be an $I$-reduced $3$-graph.  We say that
	$\A$ is \emph{$(d,\ee)$-dense} if, for every $i<j<k$ in $I$, every
	ordering $(X,Y,Z)$ of the three classes
	$\Pcal^{ij},\Pcal^{ik},\Pcal^{jk}$, and all $x\in X$ and $y\in Y$,
	\[
	\bigl|N_{\A^{ijk}}(x,y;Z)\bigr|\ge d|Z|,
	\]
	where
	\[
	N_{\A^{ijk}}(x,y;Z)
	:=
	\{z\in Z:\{x,y,z\}\in E(\A^{ijk})\}.
	\]
\end{definition}

\begin{definition}[Reduced images]\label{def:reduced-map}
	Let $F$ be a $3$-graph and let $\A$ be an $I$-reduced $3$-graph.  A
	\emph{reduced map} from $F$ to $\A$ is a pair $(\iota,\psi)$ such
	that
	\begin{enumerate}[label=\textup{(\roman*)}]
		\item $\iota:V(F)\to I$ and $\psi:\partial F\to V(\A)$;
		
		\item if $xy\in\partial F$, then $\iota(x)\ne\iota(y)$ and
		$\psi(xy)\in\Pcal^{\iota(x)\iota(y)}$;
		
		\item if $xyz\in E(F)$, then
		$\{\psi(xy),\psi(xz),\psi(yz)\}
		\in E(\A^{\iota(x)\iota(y)\iota(z)})$.
	\end{enumerate}
	If such a reduced map exists, then $\A$ contains a \emph{reduced image}
	of $F$; otherwise, $\A$ is \emph{$F$-free}.
\end{definition}

The following theorem provides the bridge
between this reduced formalism and the cherry Tur\'an density by
characterizing $\piee(F)$ in terms of $F$-free reduced $3$-graphs with
arbitrarily many indices.

\begin{theorem}[Reiher~\cite{Reiher2020}]
	\label{thm:reiher-reduced}
	Let $F$ be a $3$-graph.  Then
	\[
	\piee(F)=\sup\Bigl\{d\in[0,1]:
	\begin{array}{l}
		\text{for every }m\in\N\text{ there exists an }F\text{-free,}\\
		\text{$(d,\ee)$-dense reduced $3$-graph with an index set of size }m
	\end{array}
	\Bigr\}.
	\]
\end{theorem}

\subsection{Palette lower bounds}\label{sec:palettes}

Palette constructions go back to R\"odl's random pair-coloring
idea~\cite{Rodl1986}, while the minimum-codegree palette density for
the cherry setting was introduced by
Lamaison~\cite[Section~7]{Lamaison2024}.  In this subsection, we use the
reduced formalism above to relate this palette density to the cherry
Tur\'an density.

The following lemma shows that palette obstructions give
cherry-density lower bounds. We prove it
directly by converting any palette not admitted by $F$ into an
$F$-free reduced $3$-graph whose reduced $\ee$-density is at least the
minimum codegree of the palette.

\begin{lemma}\label{prop:palette-lower-bound}
For every $3$-graph $F$,
$\pieepal(F)\le \piee(F)$.
\end{lemma}

\begin{proof}
	Let $\mathscr P=(\mathscr C,\mathscr A)$ be a palette not admitted by
	$F$, and put $d:=\deltacd(\mathscr P)$.  We show that $d\le\piee(F)$.
	By \cref{thm:reiher-reduced}, it suffices to construct, for every
	$m\in\N$, an $F$-free $(d,\ee)$-dense reduced $3$-graph with $m$
	indices.
	
	Fix $m\in\N$ and take the index set to be $[m]$.  For every
	$ij\in\binom{[m]}{2}$, let
	$\Pcal^{ij}$ be a disjoint copy of $\mathscr C$, and write $(ij,a)$ for
	the copy of $a\in\mathscr C$.  For $i<j<k$, define the constituent
	$\A^{ijk}$ by
	\[
	\{(ij,a),(ik,b),(jk,c)\}\in E(\A^{ijk})
	\quad\Longleftrightarrow\quad
	(a,b,c)\in\mathscr A.
	\]
	Denote the resulting reduced $3$-graph by $\A_m(\mathscr P)$.
	
	The definition of $\deltacd(\mathscr P)$ gives the required reduced
	density.  Indeed, fixing vertices in the classes $\Pcal^{ij}$ and
	$\Pcal^{ik}$ leaves at least $d|\mathscr C|$ choices in
	$\Pcal^{jk}$ by $d_{12}$; the pairs of classes
	$(\Pcal^{ij},\Pcal^{jk})$ and $(\Pcal^{ik},\Pcal^{jk})$ are handled by
	$d_{13}$ and $d_{23}$, respectively.  Thus $\A_m(\mathscr P)$ is
	$(d,\ee)$-dense.
	
	It remains to prove that $\A_m(\mathscr P)$ is $F$-free.  Suppose to
	the contrary that $(\iota,\psi)$ is a reduced map from $F$ to
	$\A_m(\mathscr P)$.  For each $xy\in\partial F$, let
	$\varphi(xy)\in\mathscr C$ be the unique color whose copy in
	$\Pcal^{\iota(x)\iota(y)}$ is $\psi(xy)$.
	Choose a linear order $\prec$ on $V(F)$ by ordering the vertices
	increasingly according to $\iota$ and breaking ties arbitrarily.
	
	Let $uvw\in E(F)$ with $u\prec v\prec w$.  The reduced-map condition
	implies that $\iota(u),\iota(v),\iota(w)$ are pairwise distinct, and hence
	$\iota(u)<\iota(v)<\iota(w)$.  Moreover,
	\[
	\{\psi(uv),\psi(uw),\psi(vw)\}
	\in E\bigl(\A^{\iota(u)\iota(v)\iota(w)}\bigr).
	\]
	By the definition of the constituent,
	\[
	\bigl(\varphi(uv),\varphi(uw),\varphi(vw)\bigr)\in\mathscr A.
	\]
	Thus $(\prec,\varphi)$ is an admission witness for $F$ and $\mathscr P$,
	contrary to the choice of $\mathscr P$.  Therefore
	$\A_m(\mathscr P)$ is $F$-free.
	
	Since this construction works for every $m\in\N$,
	\cref{thm:reiher-reduced} gives
	$\deltacd(\mathscr P)\le\piee(F)$.  Taking the supremum over all
	palettes not admitted by $F$ proves the lemma.
\end{proof}

\section{Proofs of the zero characterization theorems}\label{sec:main-proof}
This section states the two key lemmas underlying our zero
characterizations and then proves
\cref{thm:main,thm:palette-zero}.  Their proofs are deferred to
\cref{sec:reduced-image,sec:palette-construction}, respectively.

The first key lemma gives the implication from max-layeredness to the
vanishing of~$\piee$.

\begin{lemma}\label{lem:max-layered-zero}
	If $F$ is a max-layered $3$-graph, then $\piee(F)=0$.
\end{lemma}

The second key lemma provides a positive lower bound for every non-max-layered $3$-graph.

\begin{lemma}\label{lem:nonmax-positive}
	If a $3$-graph $F$ is not max-layered, then we have $\pieepal(F)\ge 3^{-\,|V(F)|!\,3^{2|E(F)|}} > 0$.
\end{lemma}

Together with \cref{prop:palette-lower-bound}, this immediately implies
that $\piee(F)>0$ whenever $F$ is not max-layered.

\begin{proof}[Proof of \cref{thm:main}]
	If $F$ is max-layered, then
	\cref{lem:max-layered-zero} gives $\piee(F)=0$.
	
	Conversely, if $F$ is not max-layered, then
	\cref{lem:nonmax-positive} gives $\pieepal(F)>0$.  By
	\cref{prop:palette-lower-bound}, this implies $\piee(F)>0$.
	Hence $\piee(F)=0$ is impossible.  This proves the equivalence.
\end{proof}

The preceding proof combines the two key lemmas to characterize the
vanishing of~$\piee$.  We now turn to the palette analogue.  The
non-max-layered direction again follows from
\cref{lem:nonmax-positive}.  For the max-layered direction, we show
directly that $F$ admits every palette with positive minimum codegree by
constructing the required coloring.

\begin{proof}[Proof of \cref{thm:palette-zero}]
	If $F$ is not max-layered, then
	\cref{lem:nonmax-positive} gives $\pieepal(F)>0$. 
	
	Conversely, suppose $F$ is max-layered, witnessed by an order
	$\prec$ on $V(F)$ and a map
	$\lambda:\partial F\to\N$ such that
	$\B_{\lambda,\prec}$ is max-degenerate. Let
	$\mathscr P=(\mathscr C,\mathscr A)$ be a palette satisfying
	$\deltacd(\mathscr P)>0$.  By the definition, for every
	$a,b\in\mathscr C$ we have
	\begin{equation}
		\label{codegree-palette-conditon}
		d_{12}(a,b),\ d_{13}(a,b),\ d_{23}(a,b)\ge1.
	\end{equation}
	
	We may assume $\lambda(\partial F)=[L]$ without loss of generality.
	Then, we recursively construct a map $\kappa: [L] \to \mathscr C$ by assigning colors to $[L]$ in increasing order.  First, we choose $\kappa(1)\in\mathscr C$
	arbitrarily.
	Now let $2\le \ell \le L$, and suppose that
	$\kappa(1),\ldots,\kappa(\ell -1)$ have already been chosen so that
	\begin{equation}
		\label{cond:palette-admissibility}
		\bigl(\kappa(b_1),\kappa(b_2),\kappa(b_3)\bigr)\in\mathscr A
	\end{equation}
	for every edge-pattern
	$(b_1,b_2,b_3)\in\B_{\lambda,\prec}$ whose maximum entry is smaller than
	$\ell$.  Since $\B_{\lambda,\prec}$ is max-degenerate, there is at
	most one edge-pattern in $\B_{\lambda,\prec}$ whose maximum entry is
	$\ell$.
	
	If no such edge-pattern exists, choose $\kappa(\ell)\in\mathscr C$
	arbitrarily.  Otherwise, let
	$\bb_\ell=(b_{\ell,1},b_{\ell,2},b_{\ell,3}) \in\B_{\lambda,\prec}$ be the unique edge-pattern whose maximum
	entry is $\ell$.  The entry of $\bb_\ell$ equal to $\ell$ is unique, and its other
	two entries are smaller than $\ell$ and hence have already been
	colored.  We distinguish three cases.
	
	If $\bb_\ell=(\ell ,b_{\ell,2},b_{\ell,3})$, then by \eqref{codegree-palette-conditon}, we may therefore choose
	$c\in\mathscr C$ such that
	\[
	\bigl(c,\kappa(b_{\ell,2}),\kappa(b_{\ell,3})\bigr)\in\mathscr A
	\]
	and set $\kappa(\ell):=c$.
	Similarly, if $\bb_\ell =(b_{\ell,1}, \ell, b_{\ell,3})$, we choose
	$c\in\mathscr C$ such that $
	\bigl(\kappa(b_{\ell,1}),c,\kappa(b_{\ell,3})\bigr)\in\mathscr A$
	and set $\kappa(\ell):=c$. If $\bb_\ell=(b_{\ell,1}, b_{\ell,2}, \ell)$, we choose
	$c\in\mathscr C$ such that $\bigl(\kappa(b_{\ell,1}),\kappa(b_{\ell,2}),c\bigr)\in\mathscr A$
	and set $\kappa(\ell):=c$.
	Thus \eqref{cond:palette-admissibility} holds for every edge-pattern
	whose maximum entry is at most $\ell$, so the induction continues.
	
	After step $L$, the map $\kappa$ is defined on $[L]$, and
	\eqref{cond:palette-admissibility} holds for every edge-pattern in
	$\B_{\lambda,\prec}$.

	Define
	\[
	\varphi(xy):=\kappa\bigl(\lambda(xy)\bigr)
	\qquad\text{for every }xy\in\partial F.
	\]
	Let $uvw\in E(F)$ with $u\prec v\prec w$.  Since its edge-pattern
	belongs to $\B_{\lambda,\prec}$, by
	\eqref{cond:palette-admissibility} we have
	\[
	\bigl(\varphi(uv),\varphi(uw),\varphi(vw)\bigr)
	=
	\bigl(
	\kappa(\lambda(uv)),
	\kappa(\lambda(uw)),
	\kappa(\lambda(vw))
	\bigr)
	\in\mathscr A.
	\]
	Thus the order $\prec$ and the coloring $\varphi$ show that $F$
	admits $\mathscr P$.
	
	Thus, we have proved that every palette of positive minimum codegree is
	admitted by $F$.  Consequently, we have $\pieepal(F)=0$. This proves the equivalence.
\end{proof}

\section{Proof of \cref{lem:max-layered-zero}}
\label{sec:reduced-image}
To prove \cref{lem:max-layered-zero}, we choose, for each layer,
one representative in every pair-class of the reduced 3-graph.
These representatives are chosen in increasing order of the layers.
Once the representatives of the earlier layers are fixed, each
relevant triple of indices gives a set of possible choices for
the new representative. The reduced density condition makes
each of these sets large. A pair-class is involved in many triples,
however, so its representative must belong to all the corresponding
sets on the indices we retain.
  
  Before starting the proof, we recall the following selection lemma, which follows by combining
\cite[Lemmas~8--10]{1/27}; see also
\cite[Corollaries~4.2, 4.3 and~4.6]
{ReiherRodlSchachtJLMS2018}. 
Roughly speaking, the lemma states that if each triple has a positive proportion of available choices, then after restricting to a smaller set of indices, we can choose one vertex from each vertex class in a way that works for every remaining triple.

For $i<j<k$, define $p_1(i,j,k):=ij$, $p_2(i,j,k):=ik$, and $p_3(i,j,k):=jk$.
\begin{lemma}
	\label{lem:pair-class-selection}
	For every $n\in\N$ and $d > 0$, there exists an integer
	$N = N(n, d)$ with the following property.
	Let $\A$ be an $I$-reduced $3$-graph with $|I|\ge N$, and fix
	$r\in[3]$.  Suppose that for every $i<j<k$ in $I$ we are
	given a set 
	\[
	C_{ijk}\subseteq \Pcal^{p_r(i,j,k)}
	\qquad\text{with}\qquad
	|C_{ijk}|
	\ge
	d\,|\Pcal^{p_r(i,j,k)}|.
	\]
	Then there exist a set $I'\subseteq I$ with $|I'|=n$ and vertices
	$c_{ij}\in\Pcal^{ij}$ for all $ij\in\binom{I'}{2}$ such that
	\[
	c_{p_r(i,j,k)}\in C_{ijk}
	\qquad
	\text{for every }i<j<k\text{ in }I'.
	\]
\end{lemma}

\begin{proof}[Proof of \cref{lem:max-layered-zero}]
	Fix $d\in(0,1]$, and choose a max-layered witness $(\prec,\lambda)$ for $F$. Write $V(F)=\{v_1,\ldots,v_n\}$, where $v_1\prec v_2\prec\cdots\prec v_n$.
	We may assume $\lambda(\partial F)=[L]$ without loss of generality.
	For each $\ell\in[L]$, there is at most one edge-pattern in $\B_{\lambda,\prec}$ whose maximum entry is $\ell$.  If such
	an edge-pattern exists, denote it by
	$\bb_\ell=(b_{\ell,1},b_{\ell,2},b_{\ell,3})$, and let $r_\ell\in[3]$ be the unique position such that
	$b_{\ell,r_\ell}=\ell$.  Thus $b_{\ell,q}<\ell$ for every
	$q\in[3]\setminus\{r_\ell\}$. 
	
	Let $N(n,d)$ be a threshold function supplied by
	\cref{lem:pair-class-selection}. Define integers
	$m_{L+1},m_L,\ldots,m_1$ backwards by setting $m_{L+1}= n$ and
	\begin{equation}
	m_\ell
	=
	\begin{cases}
		N(m_{\ell+1},d),
		&\text{if an edge-pattern $\bb_\ell$ exists},\\
		m_{\ell+1},
		&\text{otherwise}.
	\end{cases}
	\end{equation}
	Let $\A$ be a $(d,\ee)$-dense $I$-reduced $3$-graph with
	$|I|\ge m_1$. Choose an arbitrary set $I_1\subseteq I$ of size $m_1$.
	Now, we recursively construct nested sets
	\[
	I_1 \supseteq I_2\supseteq\cdots\supseteq I_{L+1},
	\qquad {\rm with} \ 
	|I_\ell|=m_\ell,
	\]
	together with vertices
	\[
	c_{ij}^b\in\Pcal^{ij}
	\qquad
	\text{for every }b<\ell
	\text{ and }ij\in\binom{I_\ell}{2}.
	\]
	The construction maintains the following induction condition: for every
	edge-pattern
	$\bb=(b_1,b_2,b_3)\in\B_{\lambda,\prec}$ with
	$\max\{b_1,b_2,b_3\}<\ell$ and every $i<j<k$ in $I_\ell$,
	\begin{equation}
		\label{cond:layer-embedding}
		\left\{
		c_{ij}^{b_1},
		c_{ik}^{b_2},
		c_{jk}^{b_3}
		\right\}
		\in E(\A^{ijk}).
	\end{equation}

	For $\ell=1$, condition~\eqref{cond:layer-embedding} is vacuous since no edge-pattern in $\B_{\lambda,\prec}$ has maximum entry
	$1$. Thus we set $I_2=I_1$ and choose an arbitrary vertex
	$c_{ij}^1\in\Pcal^{ij}$ for every $ij\in\binom{I_2}{2}$.
	Consequently, \eqref{cond:layer-embedding} holds at the beginning of
	step $2$.
	
	For each $\ell\in\{2,\ldots,L\}$, the construction divides into two cases.	
	Suppose first that no edge-pattern in $\B_{\lambda,\prec}$ has maximum
	entry $\ell$. In this case, set $I_{\ell+1}=I_\ell$ and choose an arbitrary vertex
	$c_{ij}^\ell\in\Pcal^{ij}$ for every
	$ij\in\binom{I_{\ell+1}}{2}$. Since no new edge-pattern has to be
	realized at this step, condition \eqref{cond:layer-embedding} holds with $\ell + 1$.
	
	Now suppose that $\bb_\ell$ exists.  For each $q\in[3]\setminus\{r_\ell\}$, the entry $b_{\ell,q}$ is smaller than
	$\ell$, so for every $ij\in\binom{I_{\ell}}{2}$, the vertex
	$c_{ij}^{b_{\ell,q}}$ has already been chosen. 
	For every $i<j<k$ in $I_\ell$, define
	\begin{equation}
		\label{eq:available-layer-vertices}
		C_{ijk}
		:=
		\left\{
		x\in\Pcal^{p_{r_\ell}(i,j,k)}:
		\{x\}\cup
		\left\{
		c_{p_q(i,j,k)}^{b_{\ell,q}}:
		q\in[3]\setminus\{r_\ell\}
		\right\}
		\in E(\A^{ijk})
		\right\}.
	\end{equation}
	Since $\A$ is
	$(d,\ee)$-dense, by \cref{def:reduced-ee-density} we have
	\[
	|C_{ijk}|
	\ge
	d\left|\Pcal^{p_{r_\ell}(i,j,k)}\right|
	\qquad
	\text{for every }i<j<k\text{ in }I_\ell.
	\]
	
	As $|I_\ell|=N(m_{\ell+1},d)$, we may apply
	\cref{lem:pair-class-selection} to the restriction of $\A$ to
	$I_\ell$, with $m_{\ell+1}$ and $r_\ell$. We obtain a set
	$I_{\ell+1}\subseteq I_\ell$ of size $m_{\ell+1}$ and vertices
	$c_{ij}^\ell \in\Pcal^{ij}$ for all
	$ij\in\binom{I_{\ell+1}}{2}$ such that
	\[
	c_{p_{r_\ell}(i,j,k)}^\ell
	\in C_{ijk}
	\qquad
	\text{for every }i<j<k\text{ in }I_{\ell+1}.
	\]
	Then by \eqref{eq:available-layer-vertices}, the edge-pattern $\bb_\ell$ is now realized on every
	triple of $I_{\ell+1}$, thus the condition~\eqref{cond:layer-embedding} remains valid with $\ell + 1$.
	
	After the final step, $|I_{L+1}|=|V(F)|$, and
	condition~\eqref{cond:layer-embedding} applies to every edge-pattern in
	$\B_{\lambda,\prec}$. Write $I_{L+1}=\{\alpha_1<\cdots<\alpha_n\}$, define $\iota:V(F)\to I$ by
	\[
	\iota(v_i)=\alpha_i
	\qquad (i\in[n]),
	\]
	and, for every $v_iv_j\in\partial F$ with $i<j$, define
	\[
	\psi(v_iv_j)
	:=
	c_{\alpha_i\alpha_j}^{\lambda(v_iv_j)}
	\in\Pcal^{\alpha_i\alpha_j}.
	\]
	Thus, for every shadow pair $v_iv_j$ with $i<j$, we have
	\[
	\iota(v_i)\ne\iota(v_j)
	\qquad\text{and}\qquad
	\psi(v_iv_j)
	\in\Pcal^{\iota(v_i)\iota(v_j)}.
	\]
	Let $v_iv_jv_k\in E(F)$ with $i<j<k$. Its edge-pattern is $\bigl(
	\lambda(v_iv_j),
	\lambda(v_iv_k),
	\lambda(v_jv_k)
	\bigr)
	\in\B_{\lambda,\prec}$.
	Thus, applying the condition \eqref{cond:layer-embedding} to
	$\alpha_i<\alpha_j<\alpha_k$ gives
	\[
	\left\{
	\psi(v_iv_j),
	\psi(v_iv_k),
	\psi(v_jv_k)
	\right\}
	\in E(\A^{\alpha_i\alpha_j\alpha_k}).
	\]
	Hence $(\iota,\psi)$ is a reduced map from $F$ to $\A$.
	
	We have therefore shown that every $(d,\ee)$-dense reduced $3$-graph
	with at least $m_1$ indices contains a reduced image of $F$. Since $d>0$ is arbitrary,
	\cref{thm:reiher-reduced} yields $\piee(F)=0$.
\end{proof}

\section{Proof of \cref{lem:nonmax-positive}}
\label{sec:palette-construction}
In this section, we prove \cref{lem:nonmax-positive}. 
We will define the palette using three functions,
one for each position of the missing color. A triple is allowed
whenever one of its colors is the prescribed function of the
other two. On a finite color set, this guarantees positive
minimum codegree.

To exclude $F$, we first fix a vertex order and choose, for each
edge, which of the three prescriptions it is to satisfy.
These choices give a system of equations in the colors assigned
to the pairs of $F$. The type families introduced below record
which variable is determined by which two others in each equation.
We first construct functions for which each such system has no
solution, and then combine the constructions to handle all vertex
orders and choices simultaneously.

In the following, we begin by
introducing the definition of type families and explaining their connection with max-layeredness, thereby preparing the ground for the palette construction
that follows.  We then construct the desired palette, and conclude the
section by proving the auxiliary lemma on which the construction relies.

\subsection{Type families and max-layeredness}
\label{subsec:type-family}

Let $X$ be a finite set, and put $\mathcal T(X):=X^3\times[3]$.  An ordered $4$-tuple $\boldsymbol{\tau}=(h,s,t,i)\in\mathcal T(X)$ is called a \emph{type} on $X$; we call $h$ its \emph{head}, $(s,t)$ its \emph{ordered pair of tails}, and $i$ its \emph{position}.  A \emph{type family} on $X$ is a subset $\mathcal T\subseteq\mathcal T(X)$.

\begin{definition}
	\label{def:layer-realizable}
	A type family $\mathcal T\subseteq\mathcal T(X)$ is \emph{layer-realizable} if there exists a map $f:X\to\N$ such that the following conditions hold.
	\begin{enumerate}[label=(P\arabic*)]
		\item\label{type-families-item:P1}
		For every type $\boldsymbol{\tau}=(h,s,t,i)\in\mathcal T$, we have $f(h)>f(s)$ and $f(h)>f(t)$.
		\item\label{type-families-item:P2}
		For any two types $\boldsymbol{\tau}=(h,s,t,i)$ and $\boldsymbol{\tau}'=(h',s',t',i')$ in $\mathcal T$, if $f(h)=f(h')$, then $f(s)=f(s')$, $f(t)=f(t')$, and $i=i'$.
	\end{enumerate}
	Such a map $f$ is called a \emph{layer realization} of $\mathcal T$.
\end{definition}

We now associate the type family with a $3$-graph $F$. 
Let $\prec$ be an order on $V(F)$. For every edge $e=uvw\in E(F)$ with $u\prec v\prec w$, write $p_1(e):=uv$, $p_2(e):=uw$, and $p_3(e):=vw$. A \emph{branch choice function} for $(F,\prec)$ is a map $\theta:E(F)\to[3]$.

For each $e\in E(F)$, let $h_e:=p_{\theta(e)}(e)$ be the shadow pair in position $\theta(e)$, and let $(s_e,t_e)$ be the ordered pair obtained by deleting that shadow pair from the ordered triple $(p_1(e),p_2(e),p_3(e))$.  Thus
\begin{equation}
	\label{(s,t)-pair}
	(s_e,t_e):=
	\begin{cases}
		(p_2(e),p_3(e)),&\text{if }\theta(e)=1,\\
		(p_1(e),p_3(e)),&\text{if }\theta(e)=2,\\
		(p_1(e),p_2(e)),&\text{if }\theta(e)=3.
	\end{cases} 
\end{equation}
We define the type family of $F$ associated with $\prec$ and $\theta$ by
\begin{equation}
	\mathcal T_{\prec,\theta}(F) :=
	\bigl\{(h_e,s_e,t_e,\theta(e)):e\in E(F)\bigr\}
	\subseteq\mathcal T(\partial F).
\end{equation}

\begin{lemma}
	\label{lem:max-layered-type-family}
	A $3$-graph $F$ is max-layered if and only if there exist an order $\prec$ on $V(F)$ and a branch choice function $\theta:E(F)\to[3]$ such that the associated type family $\mathcal T_{\prec,\theta}(F)$ is layer-realizable.
\end{lemma}

\begin{proof}
	Suppose first that $F$ is max-layered, witnessed by an order
	$\prec$ on $V(F)$ and a map
	$\lambda:\partial F\to\N$ such that
	$\B_{\lambda,\prec}$ is max-degenerate.

	For every edge $e\in E(F)$, let $\theta(e)$ be the position of the
	unique maximum entry of its edge-pattern.  We claim that the
	map $\lambda:\partial F\to\N$ is a layer realization of
	$\mathcal T_{\prec,\theta}(F)$.
	
	For each $e\in E(F)$, the shadow pair $h_e$ occupies the position of the
	unique maximum entry of the edge-pattern of $e$.  Hence
	$\lambda(h_e)>\lambda(s_e)$ and $\lambda(h_e)>\lambda(t_e)$, so condition \ref{type-families-item:P1} holds.
	Now let $e,e'\in E(F)$ satisfy
	$\lambda(h_e)=\lambda(h_{e'})$.  The corresponding edge-patterns have
	the same maximum entry.  Since distinct elements of
	$\B_{\lambda,\prec}$ have distinct maximum entries, these two
	edge-patterns must be equal.  It follows that
	$\theta(e)=\theta(e')$, $\lambda(s_e)=\lambda(s_{e'})$, and
	$\lambda(t_e)=\lambda(t_{e'})$.
	Thus condition~\ref{type-families-item:P2} holds, and
	$\mathcal T_{\prec,\theta}(F)$ is layer-realizable.
	
	Conversely, suppose that there exist an order $\prec$ on $V(F)$,
	a branch choice function $\theta:E(F)\to[3]$, and a layer realization
	$f:\partial F\to\N$ of $\mathcal T_{\prec,\theta}(F)$.  We show
	that $(\prec,f)$ is a max-layered witness for $F$.
	By condition~\ref{type-families-item:P1}, for every edge
	$e\in E(F)$, the entry in position $\theta(e)$ is the unique maximum
	entry of its edge-pattern.
	
	Now let $e,e'\in E(F)$, and suppose that their edge-patterns with
	respect to $(\prec,f)$ have the same maximum entry.  Since these maxima
	are $f(h_e)$ and $f(h_{e'})$, respectively, we have
	$f(h_e)=f(h_{e'})$.  Condition~\ref{type-families-item:P2} gives
	$\theta(e)=\theta(e')$, $f(s_e)=f(s_{e'})$, and
	$f(t_e)=f(t_{e'})$. Hence
	the two edge-patterns are equal. Therefore
	$\B_{f,\prec}$ is max-degenerate, and $F$ is max-layered.
\end{proof}

The above lemma shows that a $3$-graph $F$ is not max-layered if
and only if, for every order $\prec$ and every branch choice
function $\theta$, the associated type family
$\mathcal T_{\prec,\theta}(F)$ is not layer-realizable.  We next introduce color systems and define compatibility
between a type family and a color system.
A \emph{color-system} is a tuple
$\mathcal S=(\C,g_1,g_2,g_3)$, where $\C$ is a finite nonempty color set and $g_i:\C^2\to \C$ for each $i\in[3]$.

\begin{definition}
	Let $\mathcal T$ be a type family on $X$, and let $\mathcal S$ be a
	color-system.  We say that $\mathcal T$ is \emph{compatible with}
	$\mathcal S$ if there exists a map $\varphi:X\to \C$ such
	that
	\begin{equation}
		\label{eq:compatibility-equations}
		\varphi(h)=g_i\bigl(\varphi(s),\varphi(t)\bigr)
		\qquad\text{for every }(h,s,t,i)\in\mathcal T.
	\end{equation}
	Such a map $\varphi$ is called a \emph{compatible color function} for
	$\mathcal T$ and $\mathcal S$.
\end{definition}

For fixed $\mathcal T$ and $\mathcal S$, the equations in
\eqref{eq:compatibility-equations} form a system of compatibility equations
whose unknown is the map $\varphi:X\to \C$.  The next lemma connects the
solvability of these compatibility equations with layer realizability and
plays a key role in the proof of
\cref{lem:nonmax-positive}.  We defer its proof to
\cref{subsec:functional-obstruction}.

\begin{lemma}
	\label{lem:functional-obstruction}
	Let $X$ be a finite set and let
	$\mathcal T\subseteq\mathcal T(X)$ be a type family.  If $\mathcal T$
	is not layer-realizable, then there exists a color-system
	$\mathcal S=(\C,g_1,g_2,g_3)$ such that
	$\mathcal T$ is not compatible with $\mathcal S$.  Moreover,
	$\C$ may be chosen so that $|\C|\le 3^{3^{|\mathcal T|}}$.
\end{lemma}

\subsection{The palette construction}
\label{subsec:palette-construction}
With \cref{lem:functional-obstruction} in hand, we now combine the resulting non-solvability into the palette construction.  

\begin{proof}[Proof of \cref{lem:nonmax-positive}]
	Suppose that $F$ is not max-layered.  Let
	\[
	\Omega
	:=
	\bigl\{
	(\prec,\theta):
	\prec\text{ is an order on }V(F)
	\text{ and }\theta\text{ is a branch choice function for }(F,\prec)
	\bigr\}.
	\]
	There are $|V(F)|!$ choices for $\prec$ and, for each such order,
	$3^{|E(F)|}$ choices for $\theta$.  Hence \[|\Omega|=|V(F)|!\,3^{|E(F)|}.\]
	
	For every $\omega=(\prec,\theta)\in\Omega$, put
	$\mathcal T_\omega:=\mathcal T_{\prec,\theta}(F)$, then
	$|\mathcal T_\omega|=|E(F)|$.  By
	\cref{lem:max-layered-type-family}, the type family $\mathcal T_\omega$
	is not layer-realizable.  Hence \cref{lem:functional-obstruction}
	provides a color-system $\mathcal S_\omega
	=
	\bigl(\C_\omega,g_1^\omega,g_2^\omega,g_3^\omega\bigr)$
	such that $\mathcal T_\omega$ is not compatible with $\mathcal S_\omega$ and $|\C_\omega|\le 3^{3^{|E(F)|}}$.
	
	Let
	\[
	\C:=\prod_{\omega\in\Omega}\C_\omega.
	\]
	Then the set $\C$ is finite and nonempty, and
	\[
	|\C|
	=\prod_{\omega\in\Omega}|\C_\omega|
	\le
	\bigl(3^{3^{|E(F)|}}\bigr)^{|V(F)|!\,3^{|E(F)|}}
	=
	3^{|V(F)|!\,3^{2|E(F)|}}.
	\]
	For
	$\boldsymbol{a}=(a_\omega)_{\omega\in\Omega}$ and
	$\boldsymbol{b}=(b_\omega)_{\omega\in\Omega}$ in $\C$, define
	$g_i:\C^2\to \C$ coordinatewise by
	\[
	g_i(\boldsymbol{a},\boldsymbol{b})
	:=
	\bigl(g_i^\omega(a_\omega,b_\omega)\bigr)_{\omega\in\Omega}
	\qquad (i\in[3]),
	\]
	and put $\mathcal S:=(\C,g_1,g_2,g_3)$. Then we
	define $\mathscr A^{(1)}:=
	\bigl\{(g_1(\boldsymbol{b},\boldsymbol{c}),\boldsymbol{b},\boldsymbol{c}):\boldsymbol{b},\boldsymbol{c}\in \C\bigr\}$, $\mathscr A^{(2)}:=
	\bigl\{(\boldsymbol{a},g_2(\boldsymbol{a},\boldsymbol{c}),\boldsymbol{c}):\boldsymbol{a},\boldsymbol{c}\in \C\bigr\}$ and $\mathscr A^{(3)}:=
	\bigl\{(\boldsymbol{a},\boldsymbol{b},g_3(\boldsymbol{a},\boldsymbol{b})):\boldsymbol{a},\boldsymbol{b}\in \C\bigr\}$, and let
	\[
	\mathscr A
	:=
	\mathscr A^{(1)}\cup\mathscr A^{(2)}\cup\mathscr A^{(3)},
	\qquad
	\mathscr P:=(\C,\mathscr A).
	\]
	Thus for all $\boldsymbol{a},\boldsymbol{b},\boldsymbol{c}\in \C$, the triples
	\[
	(\boldsymbol{a},\boldsymbol{b},g_3(\boldsymbol{a},\boldsymbol{b})),\qquad
	(\boldsymbol{a},g_2(\boldsymbol{a},\boldsymbol{c}),\boldsymbol{c}),\qquad
	(g_1(\boldsymbol{b},\boldsymbol{c}),\boldsymbol{b},\boldsymbol{c})
	\]
	belong to $\mathscr A$.  Hence
	$d_{12}(\boldsymbol{a},\boldsymbol{b})\ge1$, $d_{13}(\boldsymbol{a},\boldsymbol{c})\ge1$, and $d_{23}(\boldsymbol{b},\boldsymbol{c})\ge1$, and
	therefore
	\[
	\deltacd(\mathscr P)
	\ge\frac1{|\C|}
	\ge 3^{-\,|V(F)|!\,3^{2|E(F)|}}.
	\]
	
	We now show that $F$ does not admit $\mathscr P$.  Suppose otherwise, and let $\prec$ be an order on $V(F)$ and
	$\varphi:\partial F\to\C$ be a map such that $(\prec,\varphi)$ is an
	admission witness for $F$ and $\mathscr P$.
	For every
	$e = uvw \in E(F)$ with $u\prec v\prec w$, choose $\theta(e)\in[3]$ such that
	\[
	\bigl(
	\varphi(uv),
	\varphi(uw),
	\varphi(vw)
	\bigr)
	\in\mathscr A^{(\theta(e))}.
	\]
	If this triple belongs to more than one of the three sets, choose one
	such index arbitrarily.
	
	By the definitions of $h_e,s_e,t_e$ and
	$\mathscr A^{(1)},\mathscr A^{(2)},\mathscr A^{(3)}$, we have
	\[
	\varphi(h_e)
	=
	g_{\theta(e)}
	\bigl(\varphi(s_e),\varphi(t_e)\bigr)
	\qquad (e\in E(F)).
	\]
	Thus $\mathcal T_{\prec,\theta}(F)$ is compatible with $\mathcal S$,
	witnessed by the compatible color function $\varphi$.
	
	Let $\omega_0=(\prec,\theta)\in\Omega$, and define
	$\varphi_{\omega_0}:\partial F\to \C_{\omega_0}$ by
	\[
	\varphi_{\omega_0}(x):=\bigl(\varphi(x)\bigr)_{\omega_0}
	\qquad (x\in\partial F).
	\]
	Taking the component indexed by $\omega_0$ in the preceding compatibility
	equations gives
	\[
	\varphi_{\omega_0}(h_e)
	=
	g^{\omega_0}_{\theta(e)}
	\bigl(
	\varphi_{\omega_0}(s_e),
	\varphi_{\omega_0}(t_e)
	\bigr)
	\qquad (e\in E(F)).
	\]
	Hence $\mathcal T_{\omega_0}$ is compatible with $\mathcal S_{\omega_0}$, contradicting the choice of $\mathcal S_{\omega_0}$.
	
	Therefore $F$ does not admit $\mathscr P$, and we obtain
	\[
	\pieepal(F)\ge\deltacd(\mathscr P)\ge 3^{-\,|V(F)|!\,3^{2|E(F)|}}>0,
	\]
	as required.
\end{proof}

\subsection{Proof of \cref{lem:functional-obstruction}}
\label{subsec:functional-obstruction}

It remains to prove \cref{lem:functional-obstruction}.  The proof is by
induction on the number of types.  The inductive reduction arises when
two types have the same head and position:
condition~\ref{type-families-item:P2} then forces their corresponding tails to receive the same layer values.  We now formalize this reduction in the following definition.

\begin{definition}[Merge map]
	\label{def:merge-map}
	Let $\mathcal T$ be a type family on $X$, and suppose that
	$\mathcal T$ contains two types
	\[
	\boldsymbol{\tau}_1=(h',s_1,t_1,j)
	\qquad\text{and}\qquad
	\boldsymbol{\tau}_2=(h',s_2,t_2,j)
	\]
	with $(s_1,t_1)\ne(s_2,t_2)$.  Let $\sim$ be the equivalence relation on
	$X$ generated by $s_1\sim s_2$ and $t_1\sim t_2$.
	
	Choose one representative from each equivalence class of $\sim$,
	and let $X'\subseteq X$ be the set of these representatives.  Define
	$\rho:X\to X'$ by mapping each element of $X$ to the representative
	of its equivalence class.  We call $\rho$ a \emph{merge map}
	associated with $\boldsymbol{\tau}_1$ and $\boldsymbol{\tau}_2$.
	
	Define the type family $\mathcal T_\rho$ on $X'$ by
	\[
	\mathcal T_\rho
	:=
	\bigl\{
	(\rho(h),\rho(s),\rho(t),i):
	(h,s,t,i)\in\mathcal T
	\bigr\}.
	\]
	We call $\mathcal T_\rho$ the
	\emph{$\rho$-merge type family} of $\mathcal T$.
\end{definition}

Since $\boldsymbol{\tau}_1$ and $\boldsymbol{\tau}_2$ are distinct types in
$\mathcal T$ but induce the same type in $\mathcal T_\rho$, we have
$|\mathcal T_\rho|\le |\mathcal T|-1$.  The next lemma records the two
properties of this merge operation needed in the induction.

\begin{lemma}
	\label{lem:merge-reduction}
	Let $\mathcal T$ be a type family on $X$, and let
	$\rho:X\to X'$ be a merge map associated with two types
	\[
	\boldsymbol{\tau}_1=(h',s_1,t_1,j)
	\qquad\text{and}\qquad
	\boldsymbol{\tau}_2=(h',s_2,t_2,j)
	\]
	in $\mathcal T$ with $(s_1,t_1)\ne(s_2,t_2)$.  Then the following hold.
	\begin{enumerate}[label=\textup{(\roman*)}]
		\item
		The type family $\mathcal T$ is layer-realizable if and only if
		$\mathcal T_\rho$ is layer-realizable.
		
		\item
		If $\overline{\mathcal S}
		=(\overline{\C},\overline g_1,\overline g_2,\overline g_3)$ is a
		color-system such that
		$\mathcal T_\rho$ is not compatible with $\overline{\mathcal S}$, then there exists a color-system
		$\mathcal S=(\C,g_1,g_2,g_3)$ such that
		$\mathcal T$ is not compatible with $\mathcal S$ and
		$|\C|=|\overline{\C}|^3$.
	\end{enumerate}
\end{lemma}

\begin{proof}
	For~\textup{(i)}, first suppose that $f':X'\to\N$ is a layer
	realization of $\mathcal T_\rho$, and put $f:=f'\circ\rho$.  For every
	$(h,s,t,i)\in\mathcal T$, the type
	$(\rho(h),\rho(s),\rho(t),i)$ belongs to $\mathcal T_\rho$.
	Condition~\ref{type-families-item:P1} for $f'$ therefore gives
	$f(h)>f(s)$ and $f(h)>f(t)$.
	
	Now let $(h,s,t,i),(\widetilde h,\widetilde s,\widetilde t,
	\widetilde i)\in\mathcal T$ satisfy $f(h)=f(\widetilde h)$.  Then
	$f'(\rho(h))=f'(\rho(\widetilde h))$, so
	condition~\ref{type-families-item:P2} applied to the corresponding
	types in $\mathcal T_\rho$ gives
	\[
	i=\widetilde i,
	\qquad
	f(s)=f(\widetilde s),
	\qquad\text{and}\qquad
	f(t)=f(\widetilde t).
	\]
	Thus $f$ is a layer realization of $\mathcal T$.
	
	Conversely, suppose that $f:X\to\N$ is a layer realization of
	$\mathcal T$.  Applying condition~\ref{type-families-item:P2} to
	$\boldsymbol{\tau}_1$ and $\boldsymbol{\tau}_2$ gives
	\[
	f(s_1)=f(s_2)
	\qquad\text{and}\qquad
	f(t_1)=f(t_2).
	\]
	Since $\sim$ is generated by $s_1\sim s_2$ and $t_1\sim t_2$, the map
	$f$ is constant on every equivalence class of $\sim$.  Hence the map
	$f':X'\to\N$ defined by
	\[
	f'(\rho(x)):=f(x)
	\qquad (x\in X)
	\]
	is well-defined.  Choosing preimages in $\mathcal T$ for types in
	$\mathcal T_\rho$, conditions~\ref{type-families-item:P1} and
	\ref{type-families-item:P2} for $f'$ follow directly from the
	corresponding conditions for $f$.  Thus $f'$ is a layer realization of
	$\mathcal T_\rho$.
	
	For~\textup{(ii)}, choose $c_\ast\in\overline{\C}$ and put
	$\C:=\overline{\C}^{\,3}$.  Write
	$\boldsymbol a=(a_0,a_1,a_2)$ and
	$\boldsymbol b=(b_0,b_1,b_2)$ for elements of $\C$.  For each
	$i\in[3]$, define
	\[
	g_i(\boldsymbol a,\boldsymbol b)
	:=
	\begin{cases}
		\bigl(\overline g_i(a_0,b_0),a_0,b_0\bigr),
		&\text{if }i=j,\\
		\bigl(\overline g_i(a_0,b_0),c_\ast,c_\ast\bigr),
		&\text{if }i\ne j.
	\end{cases}
	\]
	Let $\mathcal S:=(\C,g_1,g_2,g_3)$.  Clearly
	$|\C|=|\overline{\C}|^3$.
	
	Suppose for a contradiction that $\mathcal T$ is compatible with
	$\mathcal S$, and let $\varphi:X\to\C$ be a compatible color
	function.  Let $\varphi_0:X\to\overline{\C}$ be its first-coordinate
	projection.  Taking first coordinates in the compatibility equation for
	each $(h,s,t,i)\in\mathcal T$ gives $\varphi_0(h)
	=
	\overline g_i\bigl(\varphi_0(s),\varphi_0(t)\bigr)$.
	For the two selected types, compatibility also gives
	\[
	g_j\bigl(\varphi(s_1),\varphi(t_1)\bigr)
	=
	\varphi(h')
	=
	g_j\bigl(\varphi(s_2),\varphi(t_2)\bigr).
	\]
	Comparing the second and third coordinates yields
	$\varphi_0(s_1)=\varphi_0(s_2)$ and
	$\varphi_0(t_1)=\varphi_0(t_2)$.  Hence $\varphi_0$ is constant on
	every equivalence class of $\sim$, and therefore induces a well-defined
	map $\overline\varphi:X'\to\overline{\C}$ given by
	\[
	\overline\varphi(\rho(x)):=\varphi_0(x)
	\qquad (x\in X).
	\]
	For every $(h,s,t,i)\in\mathcal T$, we have
	\[
	\overline\varphi(\rho(h))
	=
	\overline g_i\bigl(
	\overline\varphi(\rho(s)),
	\overline\varphi(\rho(t))
	\bigr).
	\]
	By the definition of $\mathcal T_\rho$, this means that
	$\overline\varphi$ is a compatible color function for
	$\mathcal T_\rho$ and $\overline{\mathcal S}$, a contradiction.
\end{proof}

Now, we give the proof of \cref{lem:functional-obstruction}.

\begin{proof}[Proof of \cref{lem:functional-obstruction}]
	We argue by induction on $n:=|\mathcal T|$.  If $n=0$, then
	$\mathcal T$ is empty and hence layer-realizable, so there is nothing
	to prove.  Let $n\ge1$, and assume that the assertion holds for every
	type family with fewer than $n$ types.
	
	We now consider the following three cases.
	
	\smallskip
	\noindent
	\emph{Case 1: a position conflict.}
	Suppose that $\mathcal T$ contains two types
	$\boldsymbol{\tau}_1=(h',s_1,t_1,j_1)$ and
	$\boldsymbol{\tau}_2=(h',s_2,t_2,j_2)$ with $j_1\ne j_2$.
	Let $\C:=\{0,1\}$, let $g_{j_2}$ be constantly $1$, and let
	$g_q$ be constantly $0$ for every $q\in[3]\setminus\{j_2\}$.
	Put $\mathcal S:=(\C,g_1,g_2,g_3)$.  If
	$\varphi:X\to\C$ were a compatible color function, then the
	compatibility equations for $\boldsymbol{\tau}_1$ and
	$\boldsymbol{\tau}_2$ would give $\varphi(h')=0$ and
	$\varphi(h')=1$, respectively, a contradiction.  Thus
	$\mathcal T$ is not compatible with $\mathcal S$, and
	$|\C|=2\le3^{3^n}$.
	
	\smallskip
	\noindent
	\emph{Case 2: a merge reduction.}
	Suppose that Case~1 does not occur, but $\mathcal T$ contains two
	types $\boldsymbol{\tau}_1=(h',s_1,t_1,j)$ and
	$\boldsymbol{\tau}_2=(h',s_2,t_2,j)$ with
	$(s_1,t_1)\ne(s_2,t_2)$.  Let $\rho:X\to X'$ be a merge map
	associated with these two types.  Since they induce the same type in
	$\mathcal T_\rho$, we have $|\mathcal T_\rho|\le n-1$.  Moreover,
	part~\textup{(i)} of \cref{lem:merge-reduction} implies that
	$\mathcal T_\rho$ is not layer-realizable.  By the induction
	hypothesis, there is a color-system
	$\overline{\mathcal S}$ with color set $\overline{\C}$ such that
	$\mathcal T_\rho$ is not compatible with $\overline{\mathcal S}$
	and $|\overline{\C}|\le 3^{3^{|\mathcal T_\rho|}}$.
	Part~\textup{(ii)} of \cref{lem:merge-reduction} then yields a
	color-system $\mathcal S$ with color set $\C$ such that
	$\mathcal T$ is not compatible with $\mathcal S$ and
	\[
	|\C|
	=|\overline{\C}|^3
	\le \bigl(3^{3^{|\mathcal T_\rho|}}\bigr)^3
	=3^{3^{|\mathcal T_\rho|+1}}
	\le 3^{3^n}.
	\]
	
	\smallskip
	\noindent
	\emph{Case 3: a dependency cycle.}
	Suppose that neither Case~1 nor Case~2 occurs.  Thus any two types
	with the same head have the same position and the same ordered pair
	of tails.  Define the \emph{dependency digraph} $D$ on vertex set
	$X$ by adding the arcs $h\to s$ and $h\to t$ for every
	$(h,s,t,i)\in\mathcal T$; a loop is regarded as a directed cycle.
	
	If $D$ were acyclic, a reverse ordering of $D$ would give
	an injective map $f:X\to[|X|]$ such that $f(x)>f(y)$ for every arc
	$x\to y$.  Thus condition~\ref{type-families-item:P1} would hold.
	Moreover, if
	$(h,s,t,i),(\widetilde h,\widetilde s,\widetilde t,
	\widetilde i)\in\mathcal T$ satisfy $f(h)=f(\widetilde h)$, then
	injectivity gives $h=\widetilde h$.  Since neither Case~1 nor Case~2
	occurs, we have $i=\widetilde i$ and
	$(s,t)=(\widetilde s,\widetilde t)$, so
	condition~\ref{type-families-item:P2} also holds.  Hence $f$ would be
	a layer realization of $\mathcal T$, a contradiction.  Therefore
	$D$ contains a directed cycle
	\[
	h_0\longrightarrow h_1\longrightarrow\cdots
	\longrightarrow h_{\ell-1}\longrightarrow h_0,
	\]
	where $h_0,\ldots,h_{\ell-1}$ are distinct, $\ell\ge1$, and all
	subscripts are read modulo $\ell$.
	
	For each $r\in\{0,\ldots,\ell-1\}$, choose a type
	$\boldsymbol{\tau}_r=(h_r,s_r,t_r,i_r)\in\mathcal T$ witnessing the
	arc $h_r\to h_{r+1}$.  Choose $\sigma_r\in\{1,2\}$ so that
	$h_{r+1}=s_r$ if $\sigma_r=1$ and $h_{r+1}=t_r$ if
	$\sigma_r=2$; if both equalities hold, choose either value.  Since
	the heads $h_0,\ldots,h_{\ell-1}$ are distinct, these selected types
	are distinct, and hence $\ell\le n$.
	
	Set $\C:=(\mathbb Z/(\ell+1)\mathbb Z)^\ell$, with coordinates
	indexed by $0,\ldots,\ell-1$.  For
	$\boldsymbol a,\boldsymbol b\in\C$ and $q\in[3]$, define
	$g_q(\boldsymbol a,\boldsymbol b)$ coordinatewise by
	\begin{equation}
		\bigl(g_q(\boldsymbol a,\boldsymbol b)\bigr)_r
		=
		\begin{cases}
			a_{r+1}+1,
			&\text{if }q=i_r\text{ and }\sigma_r=1,\\
			b_{r+1}+1,
			&\text{if }q=i_r\text{ and }\sigma_r=2,\\
			0,
			&\text{if }q\ne i_r,
		\end{cases}
	\end{equation}
	where the coordinate subscripts are read modulo $\ell$ and all
	arithmetic is performed in $\mathbb Z/(\ell+1)\mathbb Z$.  Put
	$\mathcal S:=(\C,g_1,g_2,g_3)$.
	
	Suppose for a contradiction that $\varphi:X\to\C$ is a compatible
	color function.  Considering the compatibility equation for
	$\boldsymbol{\tau}_r$ in coordinate $r$ gives
	\[
	\bigl(\varphi(h_r)\bigr)_r
	=
	\bigl(\varphi(h_{r+1})\bigr)_{r+1}+1
	\qquad (0\le r<\ell).
	\]
	Chaining these equations around the cycle yields
	\[
	\bigl(\varphi(h_0)\bigr)_0
	=
	\bigl(\varphi(h_0)\bigr)_0+\ell
	\qquad\text{in }\mathbb Z/(\ell+1)\mathbb Z,
	\]
	which is impossible since $\ell\not\equiv0\pmod{\ell+1}$.  Therefore
	$\mathcal T$ is not compatible with $\mathcal S$.  Finally,
	since $1\le\ell\le n$,
	\[
	|\C|
	=(\ell+1)^\ell
	\le(n+1)^n
	\le3^{n^2}
	\le3^{3^n}.
	\]
	
	In all three cases, we obtain a color-system $\mathcal S$ such that $\mathcal T$ is not compatible with $\mathcal S$ and $|\C|\le 3^{3^{|\mathcal T|}}$. This completes the proof. 
\end{proof}

\section{Non-jump at zero}\label{sec:not-jump}
In this section we prove \cref{thm:zero-not-jump}. 
We first define the palette density corresponding to the vertex-pair density notion.
For a palette $\mathscr P=(\mathscr C,\mathscr A)$, a position
$r\in[3]$, and a color $a\in\mathscr C$, define
\[
D_r^{\mathscr P}(a)
:=
\bigl|\{(x_1,x_2,x_3)\in\mathscr A:x_r=a\}\bigr|.
\]
The \emph{vertex-pair minimum degree} of $\mathscr P$ is
\begin{equation}
	\label{eq:vertex-pair-palette-degree}
	\delta_{\ev}(\mathscr P)
	:=
	\frac{1}{|\mathscr C|^2}
	\min_{a\in\mathscr C}\min_{r\in[3]}D_r^{\mathscr P}(a).
\end{equation} 
For a $3$-graph $F$, set
\[
\pievpal(F)=\sup\{\delta_{\ev}(\mathscr P):F\text{ does not admit }\mathscr P\}.
\] 
We shall use Lamaison's palette theorem for the vertex-pair parameter.
\begin{theorem}[{\cite[Theorem~5.1]{Lamaison2024}}]\label{thm:lamaison-ev-palette}
	For every $3$-graph $F$, we have $\piev(F)=\pievpal(F)$.
\end{theorem}

For palettes $\mathscr P=(\mathscr C_P,\mathscr A_P)$ and
$\mathscr Q=(\mathscr C_Q,\mathscr A_Q)$, a map
$f:\mathscr C_P\to\mathscr C_Q$ is a \emph{palette homomorphism} if
\[
(f(x),f(y),f(z))\in\mathscr A_Q
\qquad\text{for every }(x,y,z)\in\mathscr A_P.
\]
We write $\mathscr P\preceq\mathscr Q$ when such a map exists.  Palette
admission is monotone under this relation: if a $3$-graph admits
$\mathscr P$ and $\mathscr P\preceq\mathscr Q$, then it also admits
$\mathscr Q$.
For a palette $\mathscr P=(\mathscr C_P,\mathscr A_P)$, put
\[
\rev(\mathscr P)
:=
\bigl(\mathscr C_P,
\{(z,y,x):(x,y,z)\in\mathscr A_P\}\bigr).
\]
Reversing the witnessing vertex order shows that a $3$-graph admits
$\mathscr P$ if and only if it admits $\rev(\mathscr P)$.  Consequently,
if some $3$-graph admits $\mathscr P$ but does not admit $\mathscr Q$, then
necessarily
$\mathscr P\not\preceq\mathscr Q$ and
$\mathscr P\not\preceq\rev(\mathscr Q)$.  The following palette separation
theorem states that these two necessary conditions are also sufficient.

\begin{theorem}[{\cite[Theorem~13]{KKLT-25}}]\label{thm:palette-separation}
	Let $\mathscr P$ and $\mathscr Q$ be palettes. 
	There exists a $3$-graph $F$ which admits $\mathscr P$ but does not admit $\mathscr Q$ if and only if $\mathscr P\not\preceq\mathscr Q$ and $\mathscr P\not\preceq\rev(\mathscr Q)$.
\end{theorem}

We now define the palette product construction used below.

\begin{definition}[Palette product]\label{def:palette-product}
	Let $I$ be a nonempty finite set, and for each
	$i\in I$ let
	$\mathscr P_i=(\mathscr C_i,\mathscr A_i)$ be a palette.  Their \emph{palette product} is the palette
	\[
	\mathscr P_\times
	:=
	\prod_{i\in I}\mathscr P_i
	=
	(\mathscr C_\times,\mathscr A_\times),
	\qquad
	\mathscr C_\times:=\prod_{i\in I}\mathscr C_i,
	\]
	where for product colors
	$\boldsymbol x=(x_i)_{i\in I}$,
	$\boldsymbol y=(y_i)_{i\in I}$ and
	$\boldsymbol z=(z_i)_{i\in I}$,
	\[
	(\boldsymbol x,\boldsymbol y,\boldsymbol z)\in\mathscr A_\times
	\quad\Longleftrightarrow\quad
	(x_i,y_i,z_i)\in\mathscr A_i
	\text{ for every }i\in I.
	\]
\end{definition}

For every $i\in I$, the coordinate projection
$\mathscr C_\times\to\mathscr C_i$ is a palette homomorphism, and hence
$\mathscr P_\times\preceq\mathscr P_i$.  Moreover, for every $r\in[3]$
and every product color $\boldsymbol a=(a_i)_{i\in I}$,
$D_r^{\mathscr P_\times}(\boldsymbol a)
=
\prod_{i\in I}D_r^{\mathscr P_i}(a_i)$.
It follows that
\begin{equation}
	\label{eq:product-minimum-degree}
	\delta_{\ev}(\mathscr P_\times)
	\ge
	\prod_{i\in I}\delta_{\ev}(\mathscr P_i).
\end{equation}

The next lemma reduces all palettes above a fixed vertex-pair minimum-degree threshold to
palettes on a bounded color set.

\begin{lemma}\label{lem:bounded-palette-reduction}
	For every $0<\eta<1$, there exists an integer $n=n(\eta)$ with the
	following property.  If
	$\mathscr R=(\mathscr C_R,\mathscr A_R)$ is a palette satisfying
	$\delta_{\ev}(\mathscr R)>\eta$, then there exists a palette
	$\mathscr S=([n],\mathscr A_S)$ such that
	$\delta_{\ev}(\mathscr S)\ge\eta/2$ and
	$\mathscr S\preceq\mathscr R$.
\end{lemma}

\begin{proof}
	Choose $n=n(\eta)\ge12$ so large that $3n\exp(-\eta^2n/32)<1$.
	Select a random map $f:[n]\to\mathscr C_R$ by choosing
	$f(1),\ldots,f(n)$ independently and uniformly from $\mathscr C_R$.
	Define the palette $\mathscr S_f=([n],\mathscr A_f)$ by
	\[
	(i,j,k)\in\mathscr A_f
	\quad\Longleftrightarrow\quad
	(f(i),f(j),f(k))\in\mathscr A_R.
	\]
	Then $f$ is a palette homomorphism from $\mathscr S_f$ to
	$\mathscr R$, and hence $\mathscr S_f\preceq\mathscr R$.  It remains to show that
	$\delta_{\ev}(\mathscr S_f)\ge\eta/2$ with positive probability.
	
	Fix $i\in[n]$.  Let $Z_i^{(1)}$ count the ordered pairs
	$(j,k)$ of distinct elements of $[n]\setminus\{i\}$ for which
	$(i,j,k)\in\mathscr A_f$.  Conditional on $f(i)=a$, every such pair is
	counted with probability $\frac{D_1^{\mathscr R}(a)}{|\mathscr C_R|^2}
	\ge\eta$.
	Hence
	\[
	\mathbb E\bigl[Z_i^{(1)}\mid f(i)=a\bigr]
	\ge
	\eta(n-1)(n-2)
	\ge
	\frac{3\eta n^2}{4},
	\]
	where the last inequality uses $n\ge12$.
	
	Still conditioning on $f(i)$, the variable $Z_i^{(1)}$ depends on the
	independent random variables $f(\ell)$, $\ell\ne i$.  Changing one of
	them affects at most $2(n-2)\le2n$ counted pairs.  Therefore,
	McDiarmid's bounded-differences inequality gives
	\[
	\Pr\bigl(Z_i^{(1)}<\eta n^2/2\mid f(i)\bigr)
	\le
	\exp\left(
	-\frac{2(\eta n^2/4)^2}{(n-1)(2n)^2}
	\right)
	\le
	\exp(-\eta^2n/32).
	\]
	The same estimate holds for the analogous variables
	$Z_i^{(2)}$ and $Z_i^{(3)}$, in which $i$ occupies the second and third
	positions, respectively.  By the union bound and the choice of $n$,
	with positive probability
	\[
	Z_i^{(r)}\ge\eta n^2/2
	\qquad
	\text{for every }i\in[n]\text{ and }r\in[3].
	\]
	Since $D_r^{\mathscr S_f}(i)\ge Z_i^{(r)}$, this event implies
	$\delta_{\ev}(\mathscr S_f)\ge\eta/2$.  Fixing such an outcome of $f$
	and taking $\mathscr S:=\mathscr S_f$ proves the lemma.
\end{proof}

Fix $0<\eta<1$, and let $n=n(\eta)$ be supplied by
\cref{lem:bounded-palette-reduction}.  Let $\mathfrak P_\eta$ be the collection of
all palettes $\mathscr S=([n],\mathscr A_S)$ satisfying
$\delta_{\ev}(\mathscr S)\ge\eta/2$.  This collection is finite and nonempty, since it contains the full palette $([n],[n]^3)$.  We define the \emph{universal palette}
\[
\mathscr U_\eta
:=
\prod_{\mathscr S\in\mathfrak P_\eta}\mathscr S.
\]
By \eqref{eq:product-minimum-degree}, we have $\delta_{\ev}(\mathscr U_\eta)
\ge
(\eta/2)^{|\mathfrak P_\eta|}
>0$.
Moreover, if $\mathscr R$ is any palette with
$\delta_{\ev}(\mathscr R)>\eta$, then
\cref{lem:bounded-palette-reduction} gives some
$\mathscr S\in\mathfrak P_\eta$ with
$\mathscr S\preceq\mathscr R$, while the corresponding coordinate
projection gives $\mathscr U_\eta\preceq\mathscr S$.  Thus
\begin{equation}
	\label{eq:universal-palette}
	\mathscr U_\eta\preceq\mathscr R
	\qquad
	\text{whenever }\delta_{\ev}(\mathscr R)>\eta.
\end{equation}

We now introduce a family of special \emph{affine palettes} with positive
minimum codegree.  For every fixed palette with positive vertex-pair
minimum degree, a suitable member of this family can be separated from it.
This is the key ingredient in the proof of \cref{thm:zero-not-jump}.

\begin{definition}[Affine palette]\label{def:affine-palette}
	For a prime $p>2$, the \emph{affine palette over $\mathbb F_p$} is
	\[
	\mathscr Q_p=(\mathscr C_p,\mathscr A_p),
	\qquad
	\mathscr C_p:=\mathbb F_p,
	\qquad
	\mathscr A_p
	:=
	\{(x,y,z)\in\mathbb F_p^3:x+y-2z=1\}.
	\]
\end{definition}

Since the coefficients $1,1,-2$ are non-zero in $\mathbb F_p$, fixing any
two coordinates determines the remaining coordinate uniquely, whereas fixing
any one coordinate leaves exactly $p$ choices for the other two.  Hence
\begin{equation}
	\label{eq:affine-palette-degrees}
	\deltacd(\mathscr Q_p)
	=
	\delta_{\ev}(\mathscr Q_p)
	=
	\frac1p.
\end{equation}

\begin{lemma}\label{lem:affine-separator}
	Let $0<\eta<1$, and let $\mathscr P=(\mathscr C,\mathscr A)$ be a palette
	with $\delta_{\ev}(\mathscr P)>0$.  Then there exists a prime
	$p>\max\{2,1/\eta\}$ such that
	\[
	\mathscr P\not\preceq\mathscr Q_p
	\qquad\text{and}\qquad
	\mathscr P\not\preceq\rev(\mathscr Q_p).
	\]
\end{lemma}

\begin{proof}
	First consider the system of
	affine equations over $\mathbb Q$
	\begin{equation}\label{eq:affine-system}
		X_a+X_b-2X_c=1
		\qquad
		((a,b,c)\in\mathscr A),
	\end{equation}
	with one unknown $X_x$ for each $x\in\mathscr C$.  We claim that this system
	has no solution over $\mathbb Q$.  Indeed, suppose that
	$X:\mathscr C\to\mathbb Q$ is a solution, and choose $c_0\in\mathscr C$
	with $X_{c_0}=\max_{x\in\mathscr C}X_x$.  Since
	$\delta_{\ev}(\mathscr P)>0$, the color $c_0$ appears as the third
	coordinate of some admissible triple.  Thus there are $a,b\in\mathscr C$
	such that $(a,b,c_0)\in\mathscr A$.  The corresponding equation gives
	\[
	X_a+X_b-2X_{c_0}=1.
	\]
	But $X_a\le X_{c_0}$ and $X_b\le X_{c_0}$, so the left-hand side is at
	most $0$, a contradiction.
	
	We now pass from $\mathbb Q$ to $\mathbb F_p$.  We use the following
	standard consequence of Gaussian elimination.  If a finite linear system
	with integer coefficients has no solution over $\mathbb Q$, then there is
	an integer linear combination of its equations whose left-hand side is
	identically $0$ and whose right-hand side is a non-zero integer.  To see
	this, perform Gaussian elimination over $\mathbb Q$ on the augmented matrix;
	since the system is inconsistent, one obtains a row
	$0\ \cdots\ 0\mid r$ with $r\in\mathbb Q^\times$.  This row is a rational
	linear combination of the original equations, and clearing denominators
	gives the desired integer linear combination.
	
	Applying this observation to \eqref{eq:affine-system}, we obtain integers
	$\gamma_{abc}$, indexed by $(a,b,c)\in\mathscr A$, such that
	\[
	\sum_{(a,b,c)\in\mathscr A}
	\gamma_{abc}(X_a+X_b-2X_c)=0
	\]
	as a formal linear expression, while
	$N_1:=\sum_{(a,b,c)\in\mathscr A}\gamma_{abc}$ is a non-zero integer.
	Therefore, if $p\nmid N_1$, then \eqref{eq:affine-system} has no solution
	over $\mathbb F_p$: otherwise, reducing the above integer linear combination
	modulo $p$ would give $0=N_1\pmod p$, a contradiction.  A map $f:\mathscr C\to\mathbb F_p$ is a palette homomorphism from
	$\mathscr P$ to $\mathscr Q_p$ if and only if the assignment
	$X_x=f(x)$ gives a solution of \eqref{eq:affine-system} over
	$\mathbb F_p$.  Hence $\mathscr P\not\preceq\mathscr Q_p$ whenever
	$p\nmid N_1$.
	
	It remains to exclude homomorphisms from $\mathscr P$ to
	$\rev(\mathscr Q_p)$.  A map $f:\mathscr C\to\mathbb F_p$ is such a
	homomorphism if and only if the assignment $X_x=f(x)$ solves the reverse
	system over $\mathbb F_p$:
	\begin{equation}\label{eq:reverse-affine-system}
		X_c+X_b-2X_a=1
		\qquad
		((a,b,c)\in\mathscr A).
	\end{equation}
	The same maximum argument shows that this system has no solution over
	$\mathbb Q$.  Indeed, if $X$ were a rational solution, choose
	$a_0\in\mathscr C$ maximizing $X$.  Since $\delta_{\ev}(\mathscr P)>0$,
	there are $b,c\in\mathscr C$ such that $(a_0,b,c)\in\mathscr A$.  Then
	$X_c+X_b-2X_{a_0}=1$, while the left-hand side is at most $0$, a
	contradiction.
	
	Applying the same Gaussian-elimination observation to
	\eqref{eq:reverse-affine-system}, we obtain a non-zero integer $N_2$ such
	that the reverse system has no solution over $\mathbb F_p$ whenever
	$p\nmid N_2$.  Hence $\mathscr P\not\preceq\rev(\mathscr Q_p)$ for all such
	primes $p$.
	
	Finally, choose a prime $p>\max\{2,1/\eta\}$ such that
	$p\nmid N_1N_2$.  This is possible because only finitely many primes divide
	the non-zero integer $N_1N_2$.  This proves the lemma.
\end{proof}

\begin{proof}[Proof of \cref{thm:zero-not-jump}]
	Put $\eta:=\min\{\eps,1/2\}$ and construct the universal palette
	$\mathscr U_\eta$ above.  Since $\delta_{\ev}(\mathscr U_\eta)>0$,
	\cref{lem:affine-separator} gives a prime
	$p>\max\{2,1/\eta\}$ such that
	\[
	\mathscr U_\eta\not\preceq\mathscr Q_p
	\qquad\text{and}\qquad
	\mathscr U_\eta\not\preceq\rev(\mathscr Q_p).
	\]
	By \eqref{eq:affine-palette-degrees},
	$0<\deltacd(\mathscr Q_p)=\delta_{\ev}(\mathscr Q_p)=1/p<\eta$.
	Therefore, \cref{thm:palette-separation} gives a $3$-graph $F_\eps$ which
	admits $\mathscr U_\eta$ but does not admit $\mathscr Q_p$.
	
	We first establish the upper bound.  Let $\mathscr R$ be any palette with
	$\delta_{\ev}(\mathscr R)>\eta$.  By \eqref{eq:universal-palette},
	$\mathscr U_\eta\preceq\mathscr R$.  Since $F_\eps$ admits
	$\mathscr U_\eta$, it also admits $\mathscr R$.  Thus every palette not
	admitted by $F_\eps$ has vertex-pair minimum degree at most $\eta$, and so
	$\pievpal(F_\eps)\le\eta$.  By \cref{thm:lamaison-ev-palette}, we have $\piev(F_\eps)\le\eta\le\eps$.
	
	On the other hand, $F_\eps$ does not admit $\mathscr Q_p$.  Hence, by
	\eqref{eq:affine-palette-degrees}, we have $\pieepal(F_\eps)
	\ge
	\deltacd(\mathscr Q_p)
	=
	\frac1p$.
	It follows from \cref{prop:palette-lower-bound} that
	$\piee(F_\eps)\ge1/p$.  Combining the two bounds with the general inequality
	$\piee(F_\eps)\le\piev(F_\eps)$ gives
	\[
	0<\frac1p
	\le
	\piee(F_\eps)
	\le
	\piev(F_\eps)
	\le
	\eta
	\le
	\eps,
	\]
	as required.
\end{proof}

\section{Concluding remarks}\label{sec:conclusion}
To conclude this paper, we first provide an explicit family of $3$-graphs that separates the cherry and vertex-pair Tur\'an densities as an application of~\cref{thm:main}, and then discuss further directions concerning related positive densities.

\subsection{An explicit family}

The following family gives an explicit separation between the cherry and
vertex-pair Tur\'an densities.

\begin{lemma}\label{lem:star-plus-separation}
	For every $r\ge4$, let
	\[
	S_r^+
	=
	\Bigl([r],
	\{1ij:2\le i<j\le r\}\cup\{234\}\Bigr).
	\]
	Then we have $\piee(S_r^+)=0$ and $\piev(S_r^+)\ge\frac12$.
\end{lemma}

\begin{proof}
	The $3$-graph $S_r^+$ is the full $3$-uniform star with center $1$,
	together with the additional edge $234$; in particular,
	$S_4^+=K_4^{(3)}$.
	
	We first construct a max-layered witness for $S_r^+$.  Let $\prec$ be
	the natural order $1\prec2\prec\cdots\prec r$.  Let
	\[
	\xi:
	\{ij:2\le i<j\le r,\ j\ge5\}
	\longrightarrow
	\{4,5,6,\ldots\}
	\]
	be injective, and define $\lambda:\partial S_r^+\to\N$ by
	\[
	\lambda(12)=2,
	\qquad
	\lambda(34)=3,
	\qquad
	\lambda(1i)=1\ \ (3\le i\le r),
	\qquad
	\lambda(23)=\lambda(24)=1,
	\]
	and
	\[
	\lambda(ij)=\xi(ij)
	\qquad
	(2\le i<j\le r,\ j\ge5).
	\]
	With respect to $(\prec,\lambda)$, the edges $123$ and $124$ have the
	same edge-pattern $(2,1,1)$, while the edges $134$ and $234$ have the
	same edge-pattern $(1,1,3)$.  Every remaining edge has the form $1ij$
	with $j\ge5$, and its edge-pattern is
	$\bigl(\lambda(1i),\lambda(1j),\xi(ij)\bigr)$.  Since
	$\xi(ij)>3$, its third entry is the unique maximum entry.  The two
	edge-patterns $(2,1,1)$ and $(1,1,3)$ have maximum entries $2$ and $3$,
	respectively, while the injectivity of $\xi$ ensures that the remaining
	edge-patterns have pairwise distinct maximum entries.  Hence the maximum entry of every
	edge-pattern in $\B_{\lambda,\prec}$ occurs in a unique position, and
	distinct edge-patterns in $\B_{\lambda,\prec}$ have distinct maximum
	entries.  Thus $\B_{\lambda,\prec}$ is max-degenerate, so
	$(\prec,\lambda)$ is a max-layered witness for $S_r^+$.  By
	\cref{thm:main}, $\piee(S_r^+)=0$.

	For the vertex-pair Tur\'an density, observe that
	$K_4^{(3)}\subseteq S_r^+$ on the vertex set $\{1,2,3,4\}$.  Since
	$F\subseteq G$ implies $\piev(F)\le\piev(G)$, the theorem of Reiher,
	R\"odl and Schacht~\cite{ReiherRodlSchachtk43} gives
	\[
	\piev(S_r^+)
	\ge
	\piev(K_4^{(3)})
	=
	\frac12,
	\]
completing the proof.
\end{proof}

For $r\ge5$, this family does not arise from the known example
$K_4^{(3)}$ through the blow-up construction.  Indeed, $S_r^+$ is not
contained in any blow-up of $K_4^{(3)}$: the shadow of every $3$-graph
contained in such a blow-up is $4$-colorable, whereas
$\partial S_r^+=K_r$ has chromatic number $r>4$.  Thus the family
$(S_r^+)_{r\ge4}$ gives a separation between the cherry and vertex-pair
parameters with unbounded shadow chromatic number.  Moreover, the equality
$\piee(S_r^+)=0$ cannot be deduced from $\piee(K_4^{(3)})=0$ by
monotonicity, since $K_4^{(3)}\subseteq S_r^+$ only yields
$\piee(K_4^{(3)})\le\piee(S_r^+)$.

\subsection{Further directions}

By \cref{thm:zero-not-jump}, for every $\eps>0$ there exists a finite
$3$-graph $F_\eps$ such that
\[
0<\piee(F_\eps)\le\piev(F_\eps)\le\eps.
\]
The proof constructs a universal palette as the product of finitely many
palettes and then separates it from a suitable affine palette of positive
minimum codegree. \cref{thm:palette-separation} consequently produces a
finite $3$-graph that admits the universal palette but does not admit the
affine palette.  Although this proves that positive cherry Tur\'an densities
occur arbitrarily close to zero, it does not yield an explicit family.  It would therefore be interesting to find a natural
sequence $(G_m)_{m\ge1}$ of small or simply described $3$-graphs such that
\[
0<\piee(G_m)\longrightarrow0.
\]
Such a sequence would give a more explicit description of the positive part
of the cherry Tur\'an density near zero.

\subsection*{AI Declaration}
All mathematical ideas, arguments, results, and proofs were developed by the authors. AI tools were only used for language polishing and proofreading in this paper.

\bibliographystyle{abbrv}
\bibliography{ref}

\end{document}